




\input amstex
\documentstyle{amsppt}

\magnification=1200
\parskip 12pt
\pagewidth{5.4in}
\pageheight{7.2in}

\baselineskip=12pt
\expandafter\redefine\csname logo\string@\endcsname{}
\NoBlackBoxes                
\NoRunningHeads
\redefine\no{\noindent}

\define\C{\Bbb C}
\define\R{\Bbb R}
\define\Z{\Bbb Z}

\redefine\P{\Bbb P}

\define\al{\alpha}
 
\define\ga{\gamma}

\define\la{\lambda}

\define\Om{\Omega}

\define\th{\theta}
\define\om{\omega}
\define\Th{\Theta}
\define\Ga{\Gamma}

\define\sub{\subseteq}  
\define\es{\emptyset}
     
\redefine\b{\partial}
\define\p{\perp}     

\def\pr{\prime} 
\define\st{\ \vert\ }   

\redefine\ll{\lq\lq}
\redefine\rr{\rq\rq\ }
\define\rrr{\rq\rq}

\redefine\sin{\operatorname {sin}}
\redefine\cos{\operatorname {cos}}
\define\SkewHerm{\operatorname {SkewHerm}}
\define\diag{\operatorname {diag}}
\define\trace{\operatorname {trace}}
\redefine\dim{\operatorname {dim}}
\define\codim{\operatorname {codim}}
\define\End{\operatorname {End}}
\define\rank{\operatorname {rank}}
\define\index{\operatorname {index}}
\define\Ker{\operatorname {Ker}}
\define\Hom{\operatorname {Hom}}
\redefine\Im{\operatorname {Im}}

\define\grkrn{Gr_k(\R^n)}
\define\grkcn{Gr_k(\C^n)}
\define\glnc{GL_n\C}
\define\frakglnc{\frak{gl}_n\C}
\define\lan{\langle}
\define\ran{\rangle}
\define\llan{\lan\lan}
\define\rran{\ran\ran}
\redefine\i{{\ssize\sqrt{-1}}\,}
\define\ddt{\frac d {dt}}
\define\TC{T^{\C}}
\define\tc{\t^{\C}}
\define\GC{G^{\C}}

\redefine\u{\bold u}
\redefine\v{\bold v}
\redefine\w{\bold w}
\redefine\t{\bold t}

\topmatter
\title Morse Theory in the 1990's
\endtitle
\author Martin A. Guest
\endauthor
\endtopmatter

\document

\head
Introduction
\endhead

Since the publication of Milnor's book \cite{Mi1}
in 1963, Morse theory has been a standard topic in 
the education of geometers and topologists.  This book established
such high standards for clarity of exposition and mathematical
influence that it has been reprinted several times,
and it is still the most popular introductory reference
for the subject.  

Morse theory is not merely a useful technique.  It embodies a
far-reaching {\it idea,} which relates analysis, topology
and (most recently) physics.  This no doubt is responsible for
the resilience of Morse theory over the past several decades:
despite the essential simplicity of the idea, it seems to re-emerge
every few years to play a crucial role in some major new mathematical
development.   The title of Bott's article
\cite{Bo4} was no doubt inspired by this resilience --- here
is a subject which appears to be completely \ll worked out\rrr, yet 
which time after time has come back to yield something unexpected.
Characteristically, a few years after the appearance of \cite{Bo4},
Morse theory is again at the forefront of mathematics, as a motivational
example of a \ll topological field theory\rrr.

This article
\footnote{Submitted for publication to a volume dedicated to Brian Steer, 
to be published by Oxford University Press.
LaTeX version and figures available separately from
http://www.comp.metro-u.ac.jp/\~{}martin}
is based loosely on four lectures given at a Graduate
School in Differential Geometry, held
at the University of Durham in September 1996.  The purpose of the
lectures was to give a topical introduction to Morse theory, for
postgraduate students in the general area of geometry.
In order to save time and yet provide a concrete focus, I omitted most
of the standard proofs, and used
Morse functions on {\it Grassmannian manifolds} as a fundamental collection
of examples to illustrate the theorems.  As well as introducing some
basic aspects (such as Schubert varieties) 
of these important manifolds, this gave
the opportunity of discussing a link between Morse
theory and  Lie theory. A second feature
was that I emphasized from the start the fundamental role
played by the {\it gradient flow lines} of a Morse function.
With the benefit of hindsight --- see \S 4 --- this is a very
natural point of view to take.  It also fits well with the
Grassmannian examples, where the flow lines are known explicitly.

In late 1997
I gave a more leisurely series of lectures for
advanced undergraduate students at Tokyo Metropolitan University,
and I took this opportunity to expand greatly
my original notes. 
As in the earlier lectures, I emphasized the Grassmannians
and the role of the  gradient flow lines,
but this time I went
\ll beyond homology groups\rr in order to illustrate
the real power of Morse theory, and to give
some idea of the developments since \cite{Bo4}. I have also tried to
give a coherent account of the \ll toric\rr point of view, the full 
significance of which is only just beginning to be appreciated.

Before getting started, a few historical comments are appropriate.
After the pioneering work of Morse, the \ll modern\rr period
of Morse theory began with Bott's work in the 1950's on the
homology and homotopy groups of compact symmetric spaces.  One of
the main achievements of the Morse-theoretic approach was the
extension of this work to the loop space of a symmetric space;
this led to the discovery of the (Bott) Periodicity Theorem 
and ultimately to
K-theory. Although the role of Morse theory in this area
was quickly taken over
by the new machinery of algebraic topology, the geometrical
nature of Bott's proof of the Periodicity Theorem still retains
great appeal.
During the 1960's, Morse theory was used most prominently to
investigate the topology of manifolds, and most prominently of
all in the work of Smale, which led to a proof of the 
Poincar\acuteaccent e Conjecture in dimensions greater than four.
Following several dormant years in the early 1970's, Morse theory returned
as a guiding force in the development of mathematical Yang-Mills
theory, in which the critical points of the
Yang-Mills functional are studied.  In the
1980's, having been pushed out of the limelight by the rapidly
developing analytic and algebraic geometrical aspects of gauge
theory, Morse theory found a dramatic new role.  This time the primary
motivation was a new approach to Morse theory due to
Witten, together with an extension of these ideas by Floer.
In this approach the gradient flow lines play the central role,
and this laid the groundwork for the \ll field-theoretic\rr
point of view pioneered by Cohen-Jones-Segal, Betz-Cohen, and
Fukaya.

The main emphasis of these notes is the Morse theory of
compact --- in particular, finite-dimensional --- manifolds.
I have tried to give at least some references for each
aspect of finite-dimensional Morse theory, but not for the
infinite-dimensional theory where the subject is much more diverse.

\no{\it Acknowledgements}

I would like to thank John Bolton and Lyndon Woodward
for their invitation to the Graduate School at
the University of Durham in 1996, and Simon Salamon for encouraging me to
contribute the lecture notes to this book. 
I also wish to acknowledge financial support from
the US National Science Foundation, as well as the 
financial support and extraordinary hospitality of 
Tokyo Metropolitan University.

My own Morse-theoretical education
began at Oxford with Brian Steer in the 1970's, and I learned a lot
from my fellow students Elias Micha, 
Andrew Pressley, Simon Salamon, Pepe Seade and Socorro
Soberon. Later on I benefited greatly from discussions
about Morse theory
with Haynes Miller and Bill Richter. Most of all, however, I
have Brian to thank for suggesting Morse theory as a
suitable direction of study.
At that time I sometimes worried that
Morse theory was \ll too easy\rr a subject for serious 
mathematical research, and only much later did I understand
the (often repeated, but usually ignored) advice that 
the simplest ideas are the best ones.
Even then, with the applications
in gauge theory still on the horizon, Brian predicted that 
in due course there
would be tremendous developments in the subject.  He was right,
and for his guidance and encouragement
I dedicate these notes to him, with gratitude.

$${}$$
\head
\S 1. Morse functions
\endhead

\no{\it Brief summary:}  In this section and the next
we give a brief explanation of
Morse theory, referring mainly to the first 40 
pages of \cite{Mi1} for proofs.  We begin by defining Morse
functions and by mentioning
several nontrivial examples. We introduce 
the flow lines and the
stable and unstable manifolds, and
give some examples to illustrate these fundamental concepts.

\subheading{1.1 A basic question}

Let us agree that it is important to study manifolds.  One way
to study a manifold might be to study all possible real-valued functions
on it. Presumably, different types of manifolds will
possess different types of functions.

Every manifold admits real-valued functions, e.g. constant functions.
But it is not immediately obvious how to write down explicit formulae
for nontrivial real-valued functions on a given manifold.
For example, consider the Grassmannians
$$
\gather
\grkrn=\{
\text{real $k$-dimensional linear subspaces of}\ \R^n\}\\
\grkcn=\{
\text{complex $k$-dimensional linear subspaces of}\ \C^n\},
\endgather
$$
which are compact manifolds of (real) dimensions 
$k(n-k)$, $2k(n-k)$ respectively.  To give a function
$f:M\to\R$, where $M=\grkrn$ or $\grkcn$, we must associate
to each $k$-plane a real number. How can we do this in a
natural and nontrivial way?

As a much easier example, consider the manifold $S^1$ consisting
of complex numbers $e^{2\pi i\th}$ of unit length (i.e. the circle).
The angle $\th$ defines a function $S^1\to \R/\Z$, which is not
quite what we want, but we can obtain a real-valued function by
using $\cos2\pi\th$. It seems reasonable to regard this as
the \ll simplest\rr kind of nontrivial real-vaued function on $S^1$.
Note that this function may be interpreted as the first coordinate
of the embedding $S^1\to\R^2$, $e^{2\pi i\th}\mapsto
(\cos 2\pi\th,\sin 2\pi\th)$. This suggests a useful source
of functions on a general manifold $M$: first embed $M$ into a
euclidean space, then take a coordinate function. (But in order
to find nice functions, one has to find a nice embedding.)

As a slight modification of the previous example, we could take the
torus: $T=S^1\times S^1$.  For a fixed choice of $(a,b)\in\R^2$,
we have a \ll coordinate function\rr 
$(e^{2\pi i x},e^{2\pi i y})\mapsto 
a\cos 2\pi x  +   b\cos 2\pi y   $.
This is a coordinate function for an embedding of $T$ in $\R^4$.
We could instead use the familiar embedding of $T$ in $\R^3$, 
and take a coordinate function there (page 1 of \cite{Mi1}). But
we shall see later that such a function is slightly 
less satisfactory, for the purposes of present-day Morse theory.

The basic question which Morse theory addresses is: {\it what
is the relation between the properties of a manifold and the
properties of its real-valued functions?} By \ll properties\rr we
mean global properties, as all manifolds of the same
dimension have the same local properties. Thus, Morse theory
aims to relate topological properties of $M$ with analytical
properties of real-functions on $M$.

\subheading{1.2 Morse functions}

First we recall a standard definition:

\proclaim{Definition 1.2.1} Let $M$ be a (smooth)\footnote{The 
word \ll smooth\rrr, i.e. \ll infinitely
differentiable\rrr, will usually be omitted in future.}
manifold, and let $f:M\to\R$ be a (smooth) function.
A point $m\in M$ is called a critical point of $f$ if
$Df_m=0$. 
\endproclaim

\no The derivative $Df_m$ at $m$ is a linear functional on the
tangent space $T_mM$; thus $m$ is critical if and only if this
derivative is the zero linear functional.

In terms of a local coordinate chart $\phi:U\to \R^n$,
where $U$ is an open neighbourhood of $m$ in $M$ and
$\phi(m)=0$, $f$ 
corresponds to the function
$$
f\circ \phi^{-1}:U\to\R,
$$
and $Df_m$ is represented by the $1\times n$ matrix
$$
\left(
\frac{\b(f\circ \phi^{-1})}{\b x_1}(0),\dots,
\frac{\b(f\circ \phi^{-1})}{\b x_n}(0)
\right).
$$

It will simplify notation if we just write $f$ instead of 
$f\circ \phi^{-1}$.  Using this convention, Taylor's theorem
may be written as
$$
f(x)-f(0)=\sum a_i x_i  + \frac12 \sum a_{ij} x_i x_j
+ \text{remainder}
$$
where $a_i=\frac{\b f}{\b x_i}|_0$ and
$a_{ij} = \frac{\b^2f}{\b x_i \b x_j}|_0$.

When $m$ is a critical point, i.e. $\sum a_i x_i=0$,
it can be shown that the matrix $(a_{ij})$ defines a
symmetric bilinear form on the tangent space $T_mM$. (This bilinear
form is called the Hessian of $f$.) Hence it is
diagonalizable, and the rank and nullity do not depend on the
choice of $\phi$.  This leads to two more definitions:

\proclaim{Definition 1.2.2} Let $m$ be a critical point of $f:M\to\R$. 
The index of $m$ is defined to be the index of $(a_{ij})$, i.e.
the number of negative eigenvalues of $(a_{ij})$.
\endproclaim

\proclaim{Definition 1.2.3}  Let $m$ be a critical point of $f:M\to\R$. 
We say that $m$ is a nondegenerate critical point if and only if 
the nullity of $(a_{ij})$, i.e.
the dimension of the $0$-eigenspace of $(a_{ij})$, is zero.
\endproclaim

Since any function on a compact manifold has critical
points (e.g. maxima and minima), we cannot get very far by considering functions
without critical points. In other words, it is unreasonable to insist that
the first term in the Taylor series be a nondegenerate linear
functional at every point.  The next most favourable condition 
is that a function has no degenerate critical points, i.e. that (at
each critical point) the quadratic term in the Taylor series 
be nondegenerate:

\proclaim{Definition 1.2.4} Let $M$ be a (smooth)
manifold, and let $f:M\to\R$ be a (smooth) function. We say that
$f$ is a Morse function if and only if every critical point of $f$
is nondegenerate.
\endproclaim

To be a Morse function is in some sense a weak condition, as it
can be shown that the space of Morse functions is dense in the
space of functions. But in another sense\footnote{This somewhat 
paradoxical situation is indicative
of a \ll good\rr definition, perhaps.}
it is a strong
condition, as Morse functions have a very special local canonical form:

\proclaim{Lemma 1.2.5 (The Morse lemma)} Let $f:M\to\R$ be a
Morse function. Then, for any $m\in M$, there exists a local
chart $\phi$ at $m$ such that
$$
f(x)-f(0)=-\sum_{i=1}^{\la}x_i^2  +  \sum_{i=\la+1}^{n} x_i^2.
$$
\endproclaim

\no Note that the remainder term has disappeared.  It follows from
this formula that the index of any local maximum point is $n$,
and the index of any local minimum point is $0$.

\proclaim{Example 1.2.6 (Height functions on the torus)} 
\rm
In the diagram
below we have two embeddings of the torus $T^2$ in $\R^3$. By taking the
$z$-coordinate function, we obtain two real-valued functions on $T^2$. The critical
points are marked with crosses.

$${}$$
$${}$$
$${}$$
$${}$$
$${}$$
$${}$$
$${}$$
$${}$$
$${}$$
$${}$$

\no It is instructive to write down explicit formulae for these functions
on $T^2$, and to verify that the critical points are all
nondegenerate, with indices $2,1,1,0$.

\endproclaim

\proclaim{Example 1.2.7 (General theory of height functions)}
\rm
Let $M$ be a compact submanifold of $\R^N$. For any
$v\in \R^N$ we may define the \ll height function\rr
$h^v$ and the \ll distance function\rr $L^v$ on $M$
by
$$
h^v(m)=\lan m,v\ran, \quad
L^v(m)=\lan m-v,m-v\ran
$$
where $\lan\ ,\ \ran$ is the standard inner product. (Note
that these functions are essentially the same if $M$ is embedded 
in the sphere $S^{N-1}\sub \R^N$.) It can be shown that,
for almost all values of $v$, the functions
$h^v$ and $L^v$ are Morse functions. The critical point
theory of these functions is intimately related to the
Riemannian geometry of $M$ (with its induced metric).
A brief treatment of $L^v$ is given in \cite{Mi1};
the general theory is developed in \cite{Pa-Te}.
\qed\endproclaim

\subheading{1.3 Geometry}

To obtain geometrical information from a Morse function $f:M\to\R$
it is useful to consider the gradient vector field $\nabla f$ of $f$. 
To define the gradient vector field
we assume that a Riemannian metric $\lan\ ,\ \ran$ has been chosen
on $M$, and we define $(\nabla f)_m$ by $\lan (\nabla f)_m, X\ran
= (Df)_m(X)$ (for all $X\in T_mM$).  For many purposes (although
not all!) the particular choice of Riemannian metric is
unimportant.

By the theorem of local existence of solutions to first-order
ordinary differential equations, there exists an integral curve
$\ga$ of the vector field $-\nabla f$ through any point of $M$.
(We introduce the minus sign because we want to consider $\ga$
as flowing \ll downwards\rrr.) It should be noted that {\it explicit formulae} 
for these integral curves are not readily available in general.

If $M$ is compact,
then the domain of any such integral curve is $\R$, and $\R$ acts
on $M$ as a group of diffeomorphisms.  We shall denote this
action by $m\mapsto t\cdot m$, for $m\in M$, $t\in \R$. 
In other words, $t\cdot m = \ga(t)$, where $\ga$ is the
solution of the o.d.e. $\ga^\prime(t)=-(\nabla f)_{\ga(t)}$
such that $\ga(0)=m$. There are only
a finite number of critical points in this case (as critical
points are isolated, by the Morse lemma). Each integral curve $\ga$
\ll begins\rr and \ll ends\rr at critical points,
i.e. $\lim_{t\to\pm\infty}\ga(t)$ are critical points.
Evidently these integral curves are
constrained by the global nature of $M$, and we shall see
that they are a very useful tool for investigating $M$. 

\proclaim{Example 1.3.1}  
\rm
Let $M=S^2$, and embed $S^2$
in $\R^3$ as the unit sphere. The $z$-coordinate function
$f:S^2\to\R$,
$f(x,y,z)=z$ is a Morse function and it 
has precisely two critical points, 
the north pole and the south pole. It 
is the restriction of
$F:\R^3\to\R$, $F(x,y,z)=z$, and $-\nabla f$ is the component of
$-\nabla F=(0,0,-1)$ which is tangential to the sphere. Hence
each integral curve of $-\nabla f$ is a line of longitude running
from the north pole to the south pole (with a certain
parametrization).
\qed\endproclaim

\proclaim{Example 1.3.2
(Height functions on the torus, continued)}
\rm
Consider the two height functions on $T^2$ from Example 1.2.6.
The integral curves are illustrated below, with respect to
the induced Riemannian metric from $\R^3$. In each case the
torus is represented by $\R^2/\Z^2$, and a fundamental domain
is shown.

$${}$$
$${}$$
$${}$$
$${}$$
$${}$$
$${}$$
$${}$$
$${}$$
$${}$$
$${}$$

\endproclaim

\proclaim{Example 1.3.3}
\rm
It seems obvious that the
least number of isolated critical points of any function on $S^2$ is
$2$ (and Example 1.3.1 provides such a function). (Proof?) But
what is least number of critical points of any function on 
the torus $T^2$?  In the case of a Morse function, we shall see later
that the answer is $4$. But if arbitrary (smooth) functions are
allowed, the answer is $3$. 
(Example? See \cite{Pt}, Example 1, page 19).)
\qed\endproclaim

\subheading{1.4 Stable and unstable manifolds}

Integral curves, or flow lines, may be assembled into
the following important objects:

\proclaim{Definition 1.4.1} Let $M$ be a
compact manifold, let $f:M\to\R$ be a Morse function, and
let $m$ be a critical point of $f$.
The stable manifold $S(m)$ of $m$ is the set of points which flow
\ll down\rr to $m$, i.e.
$$
S(m)=\{x\in M\st \lim_{t\to\infty}t\cdot x=m\}.
$$
The unstable manifold $U(m)$ of $m$ is the set of points which flow
\ll up\rr to $m$, i.e.
$$
U(m)=\{x\in M\st \lim_{t\to-\infty}t\cdot x=m\}.
$$
\endproclaim 

The next result can be proved from the Morse lemma.

\proclaim{Proposition 1.4.2} Let the index of $m$ be $\la$. Then
$S(m), U(m)$ are homeomorphic (respectively) to $\R^{n-\la},\R^{\la}$.
\qed\endproclaim

It follows that a Morse function $f$ on $M$ provides
two decompositions of $M$ into disjoint \ll cells\rrr:
$$
M=\bigcup_{m \ \text{critical}\ } S(m)
=\bigcup_{m \ \text{critical}\ } U(m).
$$
We shall refer to these as the stable manifold decomposition
and the unstable manifold decomposition associated to $f$.
Observe that the unstable manifold decomposition associated to $f$
is the same as the stable manifold decomposition associated to $-f$.

By intersecting these two decompositions, we obtain a finer one,
namely
$$
M=\bigcup_{m_1,m_2 \ \text{critical}\ } F(m_1,m_2),\quad
F(m_1,m_2)=U(m_1)\cap S(m_2).
$$
This collects together the integral curves according to origin
and destination: $F(m_1,m_2)$ consists of all integral curves
which go from $m_1$ to $m_2$.

\proclaim{Example 1.4.3}  
\rm
Consider the Morse function
of Example 1.3.1, on $M=S^2$. The stable manifold decomposition
has two pieces, i.e. one copy of each of $\R^0$ and $\R^2$. The
unstable manifold decomposition is similar. 
The intersection of these decompositions
has three pieces --- two points and one copy of $\R^2-\{0\}$.
\qed\endproclaim

\proclaim{Example 1.4.4
(Height functions on the torus, continued)}
\rm
Consider again the two height functions on $T^2$ from Examples 1.2.6
and 1.3.2. In each case there are four unstable manifolds:
a $0$-cell, two $1$-cells, and a $2$-cell. The same is true for
the stable manifolds. However, the decompositions
$T^2=\cup_{m_1,m_2 \ \text{critical}\ } F(m_1,m_2)$ are quite
different.
\qed\endproclaim

The behaviour of the stable and unstable manifolds is particularly
nice if we impose the condition that they intersect transversely:

\proclaim{Definition 1.4.5} A Morse function $f:M\to\R$ on a Riemannian
manifold $M$ is said to be a Morse-Smale function if
$U(m_1)$ is transverse to $S(m_2)$ for all critical points
$m_1,m_2$ of $f$.
\endproclaim

\no For the meaning of transversality, we refer to 
Chapter 3 of \cite{Ka} or Chapter 3 of \cite{Hr}, or other
texts on differential topology.  This concept was introduced 
into Morse theory by Smale --- see \cite{Sm}, and  
the historical discussion in \cite{Bo3}.

The transversality condition implies that 
$U(m_1) \cap S(m_2)$ is a manifold, and that
$$
\codim U(m_1) \cap S(m_2) = \codim  U(m_1) 
+ \codim S(m_2)
$$
whenever $U(m_1) \cap S(m_2)$ is nonempty.  Since
$F(m_1,m_2)=U(m_1) \cap S(m_2)$, and
$\dim U(m_i)=\la_i$, we have
$$
n - \dim F(m_1,m_2) = (n - \la_1) + \la_2
$$
and hence
$$
\dim F(m_1,m_2) = \la_1 - \la_2.
$$
In particular, if there exists a flow line from $m_1$ to $m_2$,
then we must have $\la_1 > \la_2$.

\proclaim{Example 1.4.6}
\rm
Of the two Morse functions on the torus in Example 1.2.6
(see also Example 1.3.2), only one satisfies the Morse-Smale
condition. (Which one?  Note that, by the previous paragraph,
the existence of a flow line connecting two
critical points of the same index is not possible for a Morse-Smale
function.)
\qed\endproclaim

The concepts introduced so far will help us to address the
question \ll what configurations of critical points and
flow lines are possible?\rr for a Morse (or Morse-Smale)
function on a given manifold $M$.  This is a more precise
version of our original question \ll what kind of smooth
functions are possible?\rr on $M$.

Finally, we should mention that the behaviour of a smooth function
is much less predictable without the Morse condition, i.e.
in the presence of degenerate critical points. For example,
critical points are not necessarily isolated, and flow lines do not
necessarily converge to critical points. For an example of the
latter phenomenon, see page 14 of \cite{Pa-dM}.

$${}$$
\head
\S 2. Topology
\endhead

\no{\it Brief summary:} Morse theory gives a fundamental relation 
between topology and analysis.  This is traditionally expressed
by the \ll Morse inequalities\rrr. We describe various forms of
this relation, and its generalizations.

\subheading{2.1 Topology and analysis} 

The examples in \S 1 clearly suggest that there is a relation between
the topology of $M$ and the critical point data of a function $f:M\to\R$.
We list four specific examples below, following \cite{Bo3}. In
each case, topological information predicts the existence of
critical points of $f$.

\no(i) It is an elementary fact of topology that, if $M$ is compact, 
{\it then $f$ has maximum and minimum points, and these are critical points.}

\no(ii) Assume that $f$ has a finite number $k$ of critical points
(but is not necessarily a Morse function). Then we still have
$M=\cup_{m\ \text{critical}\ } S(m)$, and the $S(m)$ are a
finite number of contractible sets (but not necessarily cells).
It can be shown that this implies
the following cohomological condition:  if $m>k$, then the product
of any $m$ cohomology classes (of positive dimension) on $M$ must be zero.
We shall not make any use of this idea, which is part of the theory
of Lyusternik-Schnirelmann category, so we refer to Lecture 2 of
\cite{Bo3} for further information. However, we note that
it gives a stronger version of the prediction of (i): {\it if $M$ 
has $i$ cohomology classes (of positive dimension) whose product
is nonzero, then any function on $M$ must have at least
$i+1$ critical points.}  

For example, consider the torus $T^2$. There are two one-dimensional
cohomology classes whose product is nonzero, so any function on
$T^2$ must have at least 3 critical points. It is not possible
to find three cohomology classes (of positive dimension) 
whose product is nonzero, so we cannot improve the estimate beyond $3$. 
This is just as well, in view of Example 1.3.3.

\no(iii) The \ll Minimax Principle\rrr --- see \cite{Bo3}.

\no(iv) Finally we come to our main example, the Morse inequalities.
We shall consider this matter in detail later on, but the basic 
fact is easily stated: if $f$ is a Morse function on
a compact manifold $M$, {\it then the number of critical points
of index $k$ is greater than or equal to the $k$-th Betti number $b_k$
of $M$.} This is a considerable
generalization of (i), but in fact it is only a hint of the
power of Morse theory, as we shall see.

\subheading{2.2 The main theorem of Morse theory} 

Recall from \S 1.4 that a Morse function $f$ gives a decomposition
$$
M=\bigcup_{m \ \text{critical}\ } U(m)
$$
of a compact manifold $M$, where $U(m)$ is homeomorphic to
$\R^{\la_m}$, $\la_m$ being the index of $m$. The main theorem
of Morse theory gives information about how these pieces fit
together:

\proclaim{Theorem 2.2.1} Let $M$ be a compact manifold, and let
$f:M\to\R$ be a Morse function on $M$. Then $M$ has the homotopy
type of a cell complex, with one cell of dimension $\la$
for each critical point of index $\la$.
\endproclaim

\no To be precise, this theorem means that $M$ is homotopy
equivalent to a topological space of the form 
$$
X_r=((D^{\la_1}\cup_{f_1}D^{\la_2})\cup_{f_2}D^{\la_3})\cup_{f_3}\dots
$$
where $0=\la_1,\la_2,\la_3,\dots,\la_r=m$ are the indices of the critical
points of
$f$ (listed in increasing order, with repetitions where necessary),
and $f_1,f_2,f_3,\dots,f_{r-1}$ are certain continuous
\ll attaching maps\rrr,  with
$$
f_i:\b D^{\la_{i+1}}\to X_{i}.
$$
To simplify the notation, we shall drop the parentheses in future and 
simply write $X_r=D^{\la_1}\cup_{f_1}D^{\la_2}\cup_{f_2}
\dots\cup_{f_{r-1}}D^{\la_r}$.

In section 2.4 we sketch the main steps in the proof
of this theorem. First, however, we give some simple examples.

\proclaim{Example 2.2.2}
\rm
The height function on $S^1$ defined by the illustrated embedding
of $S^1$ in $\R^2$ has three local minima $A,B,C$ and three local
maxima $D,E,F$.

$${}$$
$${}$$
$${}$$
$${}$$
$${}$$
$${}$$
$${}$$
$${}$$
$${}$$
$${}$$

We have $S^1\simeq \{A,B,C\}\cup_{f}D^1\cup_{g}D^1\cup_{h}D^1$.
Each of the attaching maps $f,g,h$ is an injective map from
$\{-1,1\}$ to $\{A,B,C\}$.
Contemplation of this example suggests
that one should be able to say more than \ll there is at least
one minimum point and at least one maximum point\rrr. It seems plausible,
for example, that the number of local maxima must always be {\it equal} to
the number of local minima for this manifold.  We shall soon see that this
is correct (Corollary 2.3.2).
\qed \endproclaim

\proclaim{Example 2.2.3}
\rm
The height function on $S^2$ defined by the illustrated embedding
of $S^2$ in $\R^3$ has one local minimum, one critical point of
index $1$, and two local maxima.

$${}$$
$${}$$
$${}$$
$${}$$
$${}$$
$${}$$
$${}$$
$${}$$
$${}$$
$${}$$

We have $S^2\simeq D^0\cup_{f}D^1\cup_{g}D^2\cup_{h}D^2$. Here,
$D^0\cup_{f}D^1$ is a copy of $S^1$, and the maps $g,h$ serve to attach two
hemispheres to this circle. 
\qed\endproclaim

\proclaim{Example 2.2.4
(Height functions on the torus, continued)}
\rm
The two height functions on $T^2$ in Example 1.2.6 give
two cell decompositions 
$T^2\simeq D^0\cup_{f}D^1\cup_{g}D^1\cup_{h}D^2$, However,
the attaching maps behave quite differently  in each case, as is clear from
Examples 1.3.2 and 1.4.4.
\qed\endproclaim

Theorem 2.2.1 does not tell us anything about the attaching maps,
other than the fact that they exist. This appears to be a serious
disadvantage. However, we can deduce quite a lot of information on the
topology of $M$ --- for example, the Morse inequalities, 
which we discuss next --- without knowing anything further.
In the above examples, it is intuitively obvious what
the attaching maps are.

\subheading{2.3 The Morse inequalities}

Let $f:M\to\R$ be a Morse function, on a compact manifold $M$.
The Morse inequalities say that $m_i\ge b_i$, where $m_i$ is
the number of critical points of index $i$, and $b_i$ is the $i$-th
Betti number of $M$, i.e. $b_i=\dim H_i(M)$. (We use any homology
theory with coefficients in a field.) 
In the case of $T^2$, we have $b_0=b_2=1$ and $b_1=2$, so a
Morse function on $T^2$ must have at least four critical points.
We have already given examples where this happens. 

In terms of the {\it Morse
polynomial}
$$
M(t)=\sum_{i=0}^n m_i t^i
$$
and the {\it Poincar\acuteaccent e polynomial}
$$
P(t)=\sum_{i=0}^n b_i t^i
$$
we may express the Morse inequalities symbolically  as $M(t)\ge P(t)$. Using
this convenient notation, there is a stronger form of the Morse inequalities:

\proclaim{Theorem 2.3.1} Let $f:M\to\R$ be a Morse 
function, on a compact manifold $M$. Then $M(t)-P(t)=(1+t)Q(t)$,
for some polynomial $Q(t)$ such that $Q(t)\ge 0$.
\endproclaim

\no We shall sketch the proof in the next section. Before
that, we give some simple consequences and some examples.
To begin with, we put $t=-1$ in the theorem to obtain an
expression for the Euler characteristic of $M$:

\proclaim{Corollary 2.3.2} Let $f:M\to\R$ be a Morse 
function, on a compact manifold $M$. Then $\sum_{i=0}^n (-1)^i m_i
=\sum_{i=0}^n (-1)^i b_i$.
\qed\endproclaim

\no (This is a special case of the Hopf Index Theorem on vector
fields, namely for vector fields of the particular form $\nabla f$.)

\proclaim{Corollary 2.3.3}  Let $f:M\to\R$ be a Morse 
function, on a compact manifold $M$.  If $M(t)$ contains
only even powers of $t$, then $M(t)=P(t)$.
\endproclaim

\demo{Proof} Assume that $M(t)$ contains
only even powers of $t$. Then $M(t)\ge P(t)+(1+t)Q(t)$,
so neither $P(t)$ nor $(1+t)Q(t)$ can contain odd powers
of $t$. But this is possible only if $Q(t)$ is the zero polynomial.
\qed\enddemo

\no (Various generalizations of this result can be obtained by
similar reasoning --- the basic principle is that a \ll gap\rr
in the sequence $m_0,m_1,m_2,\dots$ forces a relation between
$M(t)$ and $Q(t)$. This is the \ll lacunary principle\rr of
Morse.)

\proclaim{Example 2.3.4} 
\rm
For our height functions on the torus
$T^2$, we have $M(t)=P(t)=1+2t+t^2$.
\qed\endproclaim

\proclaim{Example 2.3.5}
\rm
Let $M=\C P^n$, i.e. $n$-dimensional complex projective space.
This may be defined as $\C^{n+1}-\{0\}/\C^*$, where $\C^*$ acts
by multiplication on each coordinate of $\C^{n+1}$. It may be
identified with the space of all complex lines in  $\C^{n+1}$.
The equivalence class of $(z_0,\dots,z_n)$ will be denoted by
the standard notation $[z_0;\dots;z_n]$. 

A Morse function $f:\C P^n\to\R$ is given on page 26 of
\cite{Mi1}. It may be defined by the formula
$$
f([z_0,\dots,z_n])=\sum_{i=0}^n c_i\vert z_i\vert^2 /
\sum_{i=0}^n\vert z_i\vert^2  
$$
where $c_0<\dots<c_n$ are fixed real numbers. The critical
points of $f$ are the coordinate axes $L_0=[1;0;\dots;0]$,
$L_1=[0;1;\dots;0]$, $\dots$, $L_n=[0;0;\dots;1]$. The 
index of $L_i$ is $2i$. Since all indices are even, Corollary
2.3.4 applies, and we deduce that the Poincar\acuteaccent e
polynomial of $\C P^n$ is $1+t^2+t^4+\dots+t^{2n}$.
\qed\endproclaim

It is already clear from these results and examples that
Morse theory works \ll both ways\rrr: (1) the homology groups
of a manifold impose conditions on the critical points of
any Morse function, and (2) the critical point data of a Morse
function sometimes permits the computation of the homology groups.

If $f:M\to\R$ is a Morse function such that $M(t)=P(t)$, we
say that $f$ is {\it perfect.} The question naturally arises:
does every compact manifold posess a perfect Morse function?
This matter is discussed in detail in \cite{Pt}; the answer,
essentially, is negative. One reason is illustrated by
the next example.

\proclaim{Example 2.3.6} 
\rm
Let $\R P^n$ be 
$n$-dimensional real projective space, with $n\ge 2$.
Consider the function
$f:\R P^n\to\R$ defined by the formula of Example 2.3.5. Is
this a Morse function? (Answer: yes.) Is it a perfect Morse function?
(Answer: no, if the coefficient field is $\R$; yes,
if the coefficient field is $\Z/2\Z$.) 
\qed\endproclaim 

\subheading{2.4 Sketch proofs of the main theorems}

Theorem 2.2.1 is a consequence of two fundamental results.  Following
\cite{Bo3} we call these Theorem A and Theorem B. They describe
how the structure of the space
$$
M^t=\{ m\in M \st f(m)\le t \}
$$
changes when $t$ changes.

\proclaim{Theorem A} If $f^{-1}([a,b])$ contains no critical
point of $f$, then $M^{b}$ is diffeomorphic to $M^{a}$.
\endproclaim

\no The integral curves of (a slight modification of) $-\nabla f$
give the required diffeomorphism.  See \cite{Mi1}, pages 12--13.

\proclaim{Theorem B} If $f^{-1}([a,b])$ contains a single critical
point of $f$, then $M^{b}$ is homotopy equivalent to $M^{a}\cup_{f}D^{\la}$,
for some $f:\b D^{\la}\to M^{a}$, where $\la$ is the index of the
critical point.
\endproclaim

\no The proof is given in \cite{Mi1}, pages 14--19 (see also the
\ll Proof by picture\rr in \cite{Bo3}, pages 339--340).
By Theorem A, it suffices to consider the situation in a small
neighbourhood of the critical point.  The example $f:\R^2\to\R$,
$f(x,y)=-x^2+y^2$, in a neighbourhood of the critical point $(0,0)$,
is instructive.

To deduce Theorem 2.2.1 from Theorems A and B, an induction argument
is needed (see \cite{Mi1} again, pages 20--24). 

There are various ways of proving (and expressing) the Morse
inequalities, but all of them are based on the fact that Theorem B
allows us to compute the relative homology groups $H_{\ast}(M^{b},M^{a})$
(in the situation of Theorem B,
$H_{\ast}(M^{b},M^{a})$ is isomorphic to the coefficient field
when $\ast=\la$, and is zero otherwise). A rather formal argument
is given in \cite{Mi1}, pages 28--31; the proof in \cite{Pt} is perhaps
more transparent. Later expositions of the Morse inequalities
and the necessary background material may be found in
Chapter 8 of \cite{Ka} and Chapter 6 of \cite{Hr}, for example.)

A more intuitive version of the proof appears in \cite{Bo3}; 
the idea is to
consider how the Morse and Poincar\acuteaccent e 
polynomials $M^{a}(t)$, $P^{a}(t)$ of $M^{a}$ change when we pass
from $a$ to $b$.  Clearly we have
$$
M^{b}(t)=M^{a}(t) + t^{\la}.
$$
For the Poincar\acuteaccent e polynomial there are two possibilities,
as the $\la$-cell $D^{\la}$ either introduces a new homology
class in dimension $\la$ or bounds a homology class of $M^{a}$
in dimension $\la-1$. (More precisely, the connecting homomorphism
$H_{\la}(M^b,M^a) \to H_{\la-1}(M^a)$ in the long exact homology
sequence is either zero or injective.)
Thus we have
$$
P^{b}(t)=P^{a}(t) + \ (t^{\la} \ \text{or}\ - t^{\la-1})\,.
$$
We are interested in the difference $M^{b}(t)-P^{b}(t)$, and
for this we have
$$
M^{b}(t)-P^{b}(t)=M^{a}(t)-P^{a}(t) + 
\ (0 \ \text{or}\ t^{\la} + t^{\la-1})\,.
$$
By induction, this gives the Morse inequalities in the form
$M^{b}(t)\ge P^{b}(t)$. It also gives immediately the stronger
result $M^{b}(t)-P^{b}(t)=(1+t)Q(t)$ with $Q(t)\ge 0$.

\subheading{2.5 Generalization: Morse-Bott functions}

Since Morse functions necessarily have isolated critical points, Morse theory
immediately disqualifies many \ll natural\rr functions.   
Constant functions provide trivial examples of this phenomenon,
but there are many nontrivial ones, such as functions
which are equivariant with respect to the action of
a Lie group --- here the orbit of a critical
point consists entirely of critical points. 
We begin with a simple example:

\proclaim{Example 2.5.1}
\rm
Consider the sphere 
$S^2=\{(x,y,z)\st x^2+y^2+z^2=1\}$. Define $f:S^2\to \R$
by $f(x,y,z)=-z^2$. It is easy to verify that the critical points
are (1) the north pole $(0,0,1)$, (2) every point of the equator
$z=0$, and (3) the south pole $(0,0,-1)$. 
It is also easy to verify that (1) and (3) are nondegenerate
critical points of index zero. To understand what is happening at the
equator, let us choose local coordinates
$$
(\sqrt{1-u^2-v^2},u,v)\mapsto (u,v)
$$
around the point $(1,0,0)$. Then, with the notational conventions
of \S 1, we have
$$
f(u,v)=-v^2, \quad f(0,0)=0.
$$
Comparing this with the form of the Morse lemma,
we see that degeneracy is indicated here by the absence
of $\pm u^2$.  Now, the $u$-direction is precisely the direction
of the equator, and along the equator $f$ is constant, so degeneracy
\ll in the $u$-direction\rr is inescapable. (If $f:M\to\R$ is a
smooth function and $N$ is a connected 
submanifold of $M$ consisting of critical
points of $f$, then $f$ is constant on $N$.) In the $v$-direction,
however, we are in the situation of the Morse lemma. The integral
curves of $-\nabla f$ (with respect to the standard
induced metric on the sphere) are easy to imagine in this
situation, and we appear to have the generalized cell decomposition 
$$
S^2\simeq (\text{north pole}\cup\text{south
pole})\cup_{g}
(\text{equator}\times[-\dfrac12,\dfrac12])
$$
where 
$$
g:\text{equator}\times\{-\frac12,\frac12\}\to
\text{north pole}\cup\text{south pole}
$$ 
is the map which sends
\ll equator $\times \{ -\frac12 \}$\rr to the south pole
and \ll equator $\times \{ \frac12 \}$\rr to the north pole.
(If we use $-f$ instead of $f$, the decomposition would involve
attaching the northern and southern hemispheres to the equator.)
\qed\endproclaim 

The generalization of Morse theory to such examples was developed
by Bott in \cite{Bo1}. This depends on the following
definition:

\proclaim{Definition 2.5.2} Let $V$ be a connected submanifold
of $M$, such that every point of $V$ is a critical point of $f:M\to\R$.
We say that $V$ is a nondegenerate critical manifold of $f$ (or
NDCM, for short) if, for every $v\in V$, 
$T_vV$ is equal to the null-space $N_v$ of
the bilinear form $({\b^2f}/{\b x_i \b x_j})$ on $T_vM$.
\endproclaim 

\no With the hypotheses of the definition,
it is obvious that $T_vV$ is contained in $N_v$; for
an NDCM they are required to be  equal.  Since $V$ is assumed to
be connected, each critical point $v\in V$ has the same index, and we
refer to this number as the index of $V$. 

\proclaim{Definition 2.5.3} A function $f:M\to\R$ is a Morse-Bott
function if every critical point of $M$ belongs to a nondegenerate critical
manifold.
\endproclaim

\no Any Morse function is a Morse-Bott function, of course.
A Morse-Bott function is a Morse function if and only
if each NDCM is a point.

Bott's generalization of the main theorem of Morse theory is based
on the fact that each stable and unstable \ll cell\rr has the
structure of a vector bundle over an NDCM (see \cite{Bo1}).
In the above example, the unstable manifold of the equator is in
fact a trivial vector bundle (of rank $1$), but in general the
bundle could be nontrivial. This vector bundle is usually called
the {\it negative bundle} of the NDCM.

\proclaim{Theorem 2.5.4} Let $M$ be a compact manifold, and let
$f:M\to\R$ be a Morse-Bott function on $M$. Then $M$ has the homotopy
type of a \ll cell-bundle complex\rrr.  Each nondegenerate critical
manifold $V$ of index $\la$ contributes a cell-bundle $D^{\la}(V)$ of rank $\la$,
i.e. a fibre bundle over $V$ with fibre $D^{\la}$.
\qed\endproclaim

This means 
$$
M\simeq
D^{\la_1}(V_1)\cup_{f_1}D^{\la_2}(V_2)\cup_{f_2}D^{\la_3}(V_3)\cup_{f_3}\dots
\cup_{f_{r-1}}D^{\la_r}(V_r)
$$
where $0=\la_1,\la_2,\dots,\la_r\le m$ are the indices of the 
nondegenerate critical manifolds $V_1,V_2,\dots,V_r$ of
$f$, and $f_1,f_2,f_3,\dots,f_{r-1}$ are certain maps, with
$$
f_i:\b D^{\la_{i+1}}(V_{i+1})\to 
D^{\la_1}(V_1)\cup_{f_1}D^{\la_2}(V_2)\cup_{f_2}\dots
\cup_{f_{i-1}}D^{\la_i}(V_i).
$$
(The boundary $\b D^{\la}(V)$ of a cell-bundle $D^{\la}(V)$ is a 
sphere-bundle, i.e. a fibre bundle over $V$ with fibre $S^{\la-1}$.)

What would be an appropriate generalization of the Morse inequalities?
Well, the natural generalization of Theorem B of section 2.4 leads
one to consider the relative homology groups $H_{\ast}(D^{\la}(V),
\b D^{\la}(V))$. By the Thom Isomorphism Theorem we have
$$
H_{\ast}(D^{\la}(V),\b D^{\la}(V))\cong
H_{\ast-\la}(V;\theta)
$$
where $\theta$ is the \ll orientation sheaf\rr of the negative bundle.
There are two commonly occurring situations where this orientation
sheaf is constant: (1) if we use homology with coefficients in the
field $\Z/2\Z$, then $\theta=\Z/2\Z$; (2) if $M$ and $V$ are complex
manifolds and the negative bundle is a complex vector bundle, then
$\theta=F$ is constant for any coefficient field $F$ (or even for integer
coefficients).  If the negative bundle is
trivial (as in Example 2.5.1), then $\theta$ is certainly constant.

For notational simplicity we shall assume from now on that we are in
the favourable situation where $\theta$ is constant, so that
$H_{\ast}(D^{\la}(V),\b D^{\la}(V))\cong
H_{\ast-\la}(V)$. Following the argument
of section 2.4, we see that the
contribution of $V$ to the Poincar\acuteaccent e polynomial is either
$t^{\la}P^{V}(t)$ or $-t^{\la-1}P^{V}(t)$, 
where $P^{V}(t)$ is the Poincar\acuteaccent
e polynomial of $V$.  Let us define the
{\it Morse-Bott polynomial} of $f$ to be
$$
MB(t)=\sum_{i=1}^r P^{V_i}(t) t^{\la_i}.
$$
Then the argument of section 2.4 gives the following analogue of Theorem 2.3.1,
which we refer to as the Morse-Bott inequalities:

\proclaim{Theorem 2.5.5} Let $f:M\to\R$ be a Morse-Bott
function, on a compact manifold $M$. Then $MB(t)-P(t)=(1+t)Q(t)$,
where $Q(t)\ge 0$.
\qed\endproclaim

In particular we have $MB(t)\ge P(t)$. The comments in section
2.3 on \ll gaps\rr apply equally to this situation: if $MB(t)$ contains
no odd power of $t$, then $MB(t)=P(t)$, i.e. $f$ is a perfect Morse-Bott
function. 

Constant functions are now (satisfyingly) incorporated into the
theory, because a constant function on $M$ has one NDCM, namely
$M$ itelf, of index zero, and we have $MB(t)=P(t)$.
Somewhat more interesting are functions of the form $\pi\circ f$,
where $\pi:E\to M$ is a fibre bundle and $f:M\to\R$ is a Morse-Bott
function. Any function of this form is a Morse-Bott function; the
NDCM's are of the form $\pi^{-1}(V)$, where $V$ is an NDCM of $f$,
and the index of $\pi^{-1}(V)$ is equal to the index of $V$.

Bott explains the following intuition behind Theorem 2.5.5:
if a Morse-Bott function is deformed slightly to give a Morse
function, each NDCM $V$ of index $\la$ breaks up into a finite set of isolated
critical points, and for each Betti number $b_i$ of $V$ there
are $b_i$ critical points of index $i+\la$.  In the case of Example
2.5.1, the equator contributes $tP^{S^1}(t)=t(1+t)$, which
is equivalent to having two isolated critical points of
indices $1,2$. The way in which the equator can break up
into two isolated critical points is illustrated in the
diagrams below:

Let us look at some further concrete examples.

\proclaim{Example 2.5.6}
\rm
In Example 2.3.5 we examined the Morse function
$$
f([z_0,\dots,z_n])=\sum_{i=0}^n c_i\vert z_i\vert^2 /
\sum_{i=0}^n\vert z_i\vert^2  
$$
on $\C P^n$, where $c_0<\dots<c_n$.  If we allow {\it arbitrary}
real numbers $c_0\le\dots\le c_n$ then the same formula defines a Morse-Bott
function. Regarding $c_0,\dots,c_n$ as the eigenvalues of a diagonal
matrix $C=\diag (c_0,\dots,c_n)$, 
we can describe the NDCM's as the submanifolds
$\P(V_1),\dots,
\P(V_k)$ of $\C P^n$, where $V_1,\dots,V_k$ are the eigenspaces
of $C$. For example, in the extreme case where $c_0=0$,
$c_1=\dots=c_n=1$, we have two NDCM's, namely 
$$
\gather
\{[1;0;\dots;0]\} = \C P^0\ \text{(an isolated critical point)}\\
\{[0;*;\dots;*]\st *\in\C\} \cong \C P^{n-1}.
\endgather
$$
Clearly the isolated point is an absolute minimum point, so it has
index zero. The other NDCM consists of absolute maxima, so it has
index $2$ (i.e. $\dim \C P^n - \dim \C P^{n-1}$). Hence the Morse-Bott
polynomial is $MB(t)=1+t^2(1+t^2+\dots+t^{2(n-1)})$. This is equal
to the Poincar\acuteaccent e polynomial of $\C P^n$, so we have a
perfect Morse-Bott function. (It is easy to see that, for any choice
of $c_0\le\dots\le c_n$, we obtain a perfect Morse-Bott function, in 
fact.)
\qed\endproclaim

\proclaim{Example 2.5.7}
\rm
Let us consider a new height function on the torus $T^2$, defined
by the following embedding of $T^2$ in $\R^3$:

$${}$$
$${}$$
$${}$$
$${}$$
$${}$$
$${}$$
$${}$$
$${}$$
$${}$$
$${}$$

\no There are two NDCM's, each being a copy of $S^1$, with indices $0,1$.
If we tilt this embedding slightly, we obtain one of our old
Morse functions on $T^2$.
\qed\endproclaim 

\subheading{2.6 Generalization: noncompactness and singularities}

For a Morse function on a noncompact manifold $M$, there are
two immediate difficulties.  First, flow lines do not necessarily
exist for all time (consider the manifold 
obtained by removing a non-critical
point from a compact manifold). Second, even the flow lines which
exist for all time do not necessarily converge to critical points 
(consider the result of removing a critical point!).

If the function $f:M\to \R$ is bounded below and proper, then
the sets $f^{-1}(-\infty,a]$ are compact, and there is essentially
no difficulty in doing Morse theory. For example, it is possible
to construct Morse functions on $\C P^\infty = \cup_{n\ge 0} \C P^n$
by extending Example 2.3.5 in an obvious way.  A more interesting
example, perhaps, is the example which was fundamental for
Morse himself, namely the space of paths on a Riemannian manifold
(see \cite{Mi1}). We shall say a little more about infinite-dimensional
examples such as these in \S 4.

Morse theory can be extended in another direction, to functions
on singular spaces. Different approaches can
be found in \cite{Go-Ma}, \cite{Ki}, and section 3.2 of \cite{Kt}.

\subheading{2.7 Generalization: equivariant Morse theory}

Let $G$ be a compact Lie group acting (smoothly) on a compact
manifold $M$.  Assume that $f:M\to \R$ is $G$-equivariant,
i.e. $f(gm)=f(m)$ for all $g\in G$, $m\in M$. Then, as we
have already pointed out, $f$ is unlikely
to be a Morse function as the $G$-orbit of any critical point
consists entirely of critical points. Morse-Bott theory is
designed to cope with this kind of situation, but
there is a further extension of Morse-Bott 
theory which applies to equivariant functions.

The basic idea (see \cite{At-Bo} and \cite{Bo3})
is that one would like to relate the Morse theory of an
equivariant Morse-Bott function $f:M\to \R$ to the
\ll Morse theory\rr of the induced map $f:M/G\to \R$.
If $G$ acts freely on $M$, then $M/G$ is a compact manifold
and $M\to M/G$ is a locally trivial bundle, so such a relation
exists by the remarks following Theorem 2.5.4. But if the action of $G$
is not free, one must find another way to \ll take account of 
the symmetry due to $G$\rrr. The method introduced by
Atiyah and Bott is to replace $M/G$ by the homotopy quotient
$$
M//G = M \times_{G} EG
$$
where $EG\to BG$ is a universal bundle for $G$.
The (contractible) space $EG$ is not a finite-dimensional
manifold (if $G$ is nontrivial), but it can be constructed
as a limit of compact manifolds $(EG)_n$,
and $f$ extends naturally to a function on each compact manifold
$ M \times_{G} (EG)_n$. The following standard example
(\cite{At-Bo}) illustrates this construction.

\proclaim{Example 2.7.1}
\rm
Let $f$ be the height function on $M=S^2$ as in Example 1.3.1,
and let $G=S^1$ act on $S^2$ by \ll rotation about the vertical
axis\rrr. Thus, the two critical points are isolated
orbits, and all other orbits are circles.  Clearly ordinary
Morse theory does not work on the quotient space
$S^2/S^1 \cong [-1,1]$, so we consider instead 
$$
S^2//S^1 = S^1 \times_{S^1} S^\infty =
\lim_{n\to\infty} S^2\times_{S^1} S^{2n+1}.
$$
Now, $S^2\times_{S^1} S^{2n+1}$ is a manifold on which
(the extension of) $f$ has two nondegenerate critical
manifolds, each homeomorphic to 
$\{ \text{point} \}\times_{S^1} S^{2n+1} \cong \C P^n$.
The indices of these NDCM's are $0$ and $2$, and so we 
have a perfect Morse-Bott function with Morse-Bott
polynomial
$$
t^0(1+t^2+\dots + t^{2n}) + t^2(1+t^2+\dots + t^{2n}).
$$
Taking the limit $n\to\infty$, we obtain the Poincar\acuteaccent e
polynomial of $S^2//S^1$ as $(1+t^2)\sum_{i\ge 0} t^{2i}$
(conveniently abbreviated as $(1+t^2)/(1-t^2)$).  This agrees with
the well known fact that $S^2//S^1$ is homotopy equivalent
to $\C P^\infty \vee \C P^\infty$.  The  Poincar\acuteaccent e
polynomial of $S^2\times_{S^1} S^{2n+1}$ could also be
obtained by using the fact that $S^2\times_{S^1} S^{2n+1}$
is a bundle over $\C P^n$ with fibre $S^2$; if we choose
any perfect Morse function on $\C P^n$ then we
obtain a perfect Morse-Bott function on $S^2\times_{S^1} S^{2n+1}$
by the remarks following Theorem 2.5.4.
\qed\endproclaim

The next instructive example is also taken from \cite{At-Bo}.

\proclaim{Example 2.7.2}
\rm
Let $f:S^2\to \R$ be the function $-z^2$ of Example 2.5.1.
This Morse-Bott function is not perfect, and neither is the
extended Morse-Bott function on $S^2\times_{S^1} S^{2n+1}$.
The NDCM's of the latter are
$$
\align
\{ \text{point} \} \times_{S^1} S^{2n+1} &\cong \C P^n
\quad (\text{index}\ 2)\\
\{ \text{equator} \} \times_{S^1} S^{2n+1} &\cong S^{2n+1}
\quad (\text{index}\ 0)\\
\{ \text{point} \} \times_{S^1} S^{2n+1} &\cong \C P^n
\quad (\text{index}\ 2)
\endalign
$$
so the Morse-Bott polynomial is
$$
t^0(1+t^{2n+1}) + 2t^2(1+t^2+\dots+t^{2n}),
$$
which has a superfluous term $t^{2n+1}$. But as $n\to\infty$,
the effect of this term disappears, and we obtain a 
perfect \ll Morse-Bott function\rr on $S^2//S^1$.
\qed\endproclaim 

In general, a $G$-equivariant Morse-Bott function
$f:M\to \R$ is converted by the above procedure into
a \ll Morse-Bott function\rr on $M//G$; each NDCM  $N$
in $M$ of index $\la$ is converted into an NDCM $N//G$
of index $\la$.  Atiyah and Bott give a criterion for the
new Morse-Bott function to be perfect (even when $f$ is
not), the \ll Self-Completion Principle\rrr.

It is important to understand what this extension of
Morse theory does {\it not} do: it does not give additional
information about the homology groups of $M$, but rather about
those of the (much larger) auxiliary space $M//G$.  
Sometimes it is useful to associate a
{\it smaller} auxiliary space to the original Morse-theoretic
data, especially when dealing with infinite-dimensional
manifolds.  One example of this is provided by the
space of \ll broken geodesics\rr in Morse's original
application to the space $M=\Om X$ of closed smooth paths on
a manifold $X$, or (analogously) the space of \ll algebraic
loops\rr when $X$ is a Lie group (see \cite{Py}).
Another famous example is due to Floer (\cite{Fl}).

In all of these situations, one may regard the homology groups 
of the auxiliary construction as generalized homology groups
of the original space $M$; from equivariant Morse theory
one obtains equivariant homology groups, and from Floer's theory
Floer homology groups.
A second lesson from equivariant Morse theory is that one can
sometimes \ll do Morse theory without a Morse function\rrr:
the space $M//G$ is not a manifold in general, yet we have
constructed a Morse-Bott polynomial and the Morse-Bott inequalities
are satisfied. 

$${}$$
\head
\S 3. Morse functions on Grassmannians
\endhead

\no{\it Brief summary:} We  examine in detail a special but
very important example, namely the Morse theory of a family of
functions on the complex Grassmannian manifold $\grkcn$.  Our main tool
is an explicit formula for the integral curves.  This will
give information on the homology and cohomology of $\grkcn$, and on the behaviour
of its Schubert subvarieties.

\subheading{3.1 Morse functions on Grassmannians}

The (complex) Grassmannian
$$
\grkcn=\{
\text{complex $k$-dimensional linear subspaces of}\ \C^n\}
$$
provides a good \ll test case\rr for the theory of the previous two
sections. In addition, $\grkcn$ is an important manifold which appears
in many parts of mathematics.  We shall study it in detail in this section.

To produce a \ll nice\rr real-valued function on $\grkcn$, we use
the embedding
$$
\align
\grkcn &\sub \SkewHerm_n
=\{n\times n \ \text{complex matrices}\ A\st A^\ast = -A \}\\
V & \mapsto \i \pi_V
\endalign
$$
where $\pi_V:\C^n\to\C^n$ denotes \ll orthogonal projection on $V$\rr
with respect to the standard Hermitian inner product on $\C^n$\rrr.
The (real) vector space $\SkewHerm_n$ has an inner product 
$\llan\ ,\ \rran$, defined by $\llan A,B\rran=\trace AB^\ast$
$=-\trace AB$.
We can obtain real-valued functions on $\grkcn$ by taking height
functions with respect to this embedding. Specifically, we shall 
choose real numbers $a_1,\dots,a_n$ and consider the function
$$
f:\grkcn\to \R, \quad
V\mapsto \llan\ \i\pi_V,\i D \rran = \trace \pi_V D
$$
where $D=\diag(a_1,\dots,a_n)$ denotes the diagonal matrix
whose diagonal entries are $a_1,\dots,a_n$.

In the special case $k=1$ (so that $\grkcn=\C P^{n-1}$), we may
write
$$
V=
\left[
\matrix
z_1\\ \vdots\\ z_n
\endmatrix
\right]=
[z]\in \C P^{n-1}
$$
and we have
$$
\pi_V=\frac{zz^\ast}{\vert z\vert^2}.
$$
So our real-valued function $f:\C P^{n-1}\to \R$ is given by
$$
f([z])=\sum_{i=0}^n a_i\vert z_i\vert^2 /
\sum_{i=0}^n\vert z_i\vert^2. 
$$
This is our old friend from Example 2.3.5.

Let $V_i=\C e_i$, i.e. the line spanned by the $i$-th basis vector
of $\C^n$. We shall prove:

\proclaim{Theorem 3.1.1} Assume that $a_1>\dots>a_n\ge 0$. Then
$f$ is a perfect Morse function. An element $V$ of $\grkcn$  is a critical
point of $f$ if and only if $V=V_{i_1}\oplus\dots\oplus V_{i_k}$ for
some $i_1,\dots,i_k$ with $1\le i_1<\dots<i_k\le n$.
\endproclaim

\no One way (and indeed the conventional way) to prove this theorem
would be by direct calculation, using a local chart. We shall take
a slightly different point of view here, by concentrating on the
integral curves of $-\nabla f$. Remarkably, there is an
{\it explicit formula} for these integral curves, which greatly simplifies
the Morse-theoretic analysis of $f$:

\proclaim{Lemma 3.1.2} The integral curve $\ga$ of  $-\nabla f$
through a point $V\in\grkcn$ is given by
$\ga(t)=\diag(e^{-a_1t},\dots,e^{-a_nt})V = e^{-tD}V$.
\endproclaim

\no{\it Remark:} 
It follows immediately from this lemma that the critical points of $f$ are
the special $k$-planes $V_{\u}=V_{u_1}\oplus\dots\oplus V_{u_k}$,
as stated in the above theorem.

\demo{Proof} It is easy to verify that the tangent space of
$\grkcn$ is given as a subspace of $\SkewHerm_n$ by
$$
T_V\grkcn = \{ T-T^\ast \st T \in \Hom(V,V^\p) \},
$$
and that the orthogonal projection of any $X\in \SkewHerm_n$
on this tangent space is $\pi_V X \pi_V^\p + \pi_V^\p X \pi_V$.
Since $f$ is the restriction of the linear function
$X\mapsto \llan X, \i D\rran$ on $\SkewHerm_n$, the gradient
$\nabla f_V$ is the projection of $\i D$ on 
$T_V\grkcn$, so\footnote{This formula shows that 
$V$ is a critical point of $f$ if and
only if $DV\sub V$, $DV^\p \sub V^\p$, i.e. if and only if $V$ is of
the form $V_{i_1}\oplus\dots\oplus V_{i_k}$}
we have
$$
-\nabla f_V = -\i(\pi_V D \pi_V^\p + \pi_V^\p D \pi_V).
$$

To find $\dot\ga(t)$, we first write $e^{-tD} = U_tP_t$,
where $U_t$ is unitary and $P_t$ is an invertible complex
$n\times n$ matrix such that $P_t V = V$ (such a factorization
may be accomplished by the Gram-Schmidt orthogonalization
process).  Using this we have
$$
\align
\dot \ga(t) &=
\frac{d}{dt} \i \pi_{e^{-tD}V} \\
&=\frac{d}{dt} \i \pi_{U_t V} \\
&=\frac{d}{dt} \i U_t \pi_V U_t^{-1}\\
&= \i \{ \dot U_t \pi_V U_t^{-1} - U_t \pi_V U_t^{-1} \dot U_t U_t^{-1} \}\\
&= \i U_t \{ U_t^{-1} \dot U_t \pi_V - \pi_V U_t^{-1} \dot U_t \} U_t^{-1}\\
&= \i U_t \{ \pi_V^\p U_t^{-1} \dot U_t \pi_V - 
\pi_V U_t^{-1} \dot U_t \pi_V^\p \} U_t^{-1}.
\endalign
$$

From the identity
$$
-D = \frac{d}{dt}(e^{-tD})(e^{-tD})^{-1} =
\frac{d}{dt}(U_t P_t) (U_t P_t)^{-1}
$$
we find that $U_t^{-1} \dot U_t = - U_t^{-1} D U_t - \dot P_t P^{-1}$.
From the fact that $P_t V=V$ we have $\pi_V^\p P_t \pi_V = 0$ and hence
$$
0 = \pi_V^\p \dot P_t \pi_V = \pi_V^\p \dot P_t P_t^{-1}\pi_V.
$$
Combining these we obtain
$$
\pi_V^\p U_t^{-1} \dot U_t \pi_V = - \pi_V^\p U_t^{-1} D U_t \pi_V,
$$
so the first term of $\dot \ga(t)$ is
$$
\i U_t ( - \pi_V^\p U_t^{-1} D U_t \pi_V ) U_t^{-1}
= -\i \pi_{U_t V}^\p D \pi_{U_t V},
$$
which agrees with the first term of $-\nabla f_{\ga(t)}$.

To deal with the second term of $\dot \ga(t)$, we use the
factorization $e^{-tD} = (e^{-tD})^\ast = P_t^\ast U_t^\ast
= Q_t U_t^{-1}$, where $Q_t = P_t^\ast$. From the identity
$$
-D = (e^{-tD})^{-1}\frac{d}{dt}(e^{-tD}) =
U_t Q_t^{-1} \frac{d}{dt}(Q_t U_t^{-1}) 
$$
we find that $U_t^{-1} \dot U_t = Q_t^{-1} \dot Q_t +
U_t^{-1} D U_t$. Since $Q_t V^\p = V^\p$, we obtain
$\pi_V Q_t^{-1} \dot Q_t \pi_V^\p = 0$.  Hence the second
term of $\dot \ga(t)$ is
$$
-\i U_t ( \pi_V U_t^{-1} D U_t \pi_V^\p ) U_t^{-1} =
-\i \pi_{U_t V} D \pi_{U_t V}^\p,
$$
and this agrees with the second term of $-\nabla f_{\ga(t)}$.
\qed\enddemo 

The lemma allows us to identify the stable and unstable manifolds.
We begin with the case $k=1$. Observe that
$$
\lim_{t\to\infty}
\pmatrix
e^{-a_1t} &  &  \\
  &  \ddots &  \\
  &  & e^{-a_nt}
\endpmatrix
\left[
\matrix
* \\ \vdots \\ * \\ 1 \\ 0 \\ \vdots \\ 0
\endmatrix
\right]=
\left[
\matrix
0 \\ \vdots \\ 0 \\ 1 \\ 0 \\ \vdots \\ 0
\endmatrix
\right]=
V_{u}.
$$
Hence the stable manifold $S_{u}$ of $V_{u}$ contains all points
of the form $[* \cdots * 1 0 \cdots 0]^t$.  As $u$ varies from $1$
to $n$, such points account for all of $\C P^{n-1}$ --- so we deduce
that 
$$
S_u\cong
\left\{
\pmatrix
* \\ \vdots \\ * \\ 1 \\ 0 \\ \vdots \\ 0
\endpmatrix
\ 
\text{with}
\ 
*\in\C
\right\}
\cong \C^{u-1}.
$$
Similarly we have 
$$
U_u\cong
\left\{
\pmatrix
0 \\ \vdots \\ 0 \\ 1 \\ * \\ \vdots \\ *
\endpmatrix
\ 
\text{with}
\ 
*\in\C
\right\}
\cong \C^{n-u}.
$$
It is clear from this description that $S_u$ and $U_u$ meet tranversally
at $V_u$, and that $f$ must be a Morse function. Finally, since
the stable and unstable manifolds are even-dimensional, $f$ is perfect.

We use similar notation in the case of general $k$.  First, we represent
the critical point $V_{\u}=V_{u_1}\oplus\dots\oplus V_{u_k}$ by an
$n\times k$ matrix:
$$
\left[
\matrix
\vdots &  &  \vdots\\
V_{u_1} & \dots & V_{u_k}\\
\vdots &  & \vdots
\endmatrix
\right]
=
\left[
\matrix
  &  &  \\
1 &  &  \\
  & \ddots & \\
  &  & 1\\
  &  &
\endmatrix
\right].
$$
Then we observe that
$$
\lim_{t\to\infty}
\pmatrix
e^{-a_1t} &  &  \\
  &  \ddots &  \\
  &  & e^{-a_nt}
\endpmatrix
\left[
\matrix
* &  & * \\
\vdots &  & \vdots\\
* &  & \vdots\\
1 & & \vdots\\
  & \ddots & *\\
  &  &  1\\
  &  &  
\endmatrix
\right]
=
\left[
\matrix
  &  &  \\
1 &  &  \\
  & \ddots & \\
  &  & 1\\
  &  &
\endmatrix
\right].
$$
But any $n\times k$ matrix may be brought into the \ll echelon form\rr
$$
\left[
\matrix
* &  & * \\
\vdots &  & \vdots\\
* &  & \vdots\\
1 & & \vdots\\
  & \ddots & *\\
  &  &  1\\
  &  &  
\endmatrix
\right]
$$
by \ll column operations\rrr, so we have found the stable manifold $S_{\u}$
of $V_{\u}$. 

Now, an element of $S_{\u}$ may be represented by more than one matrix in 
echelon form.  To obtain a unique representation, we should bring the
matrix into \ll reduced echelon form\rrr, so that it looks (for
example) like this:
$$
\pmatrix
* & * & *\\
* & * & *\\
1 & 0 & 0\\
  & 1 & 0\\
  &  & *\\
  &  & 1\\
  &  & 0\\
  &  & 0
\endpmatrix.
$$
It follows that
$$
\align
\dim_{\C} S_{\u} &=
(u_1-1)+(u_2-2)+(u_3-3)+\dots+(u_k-k)\\
&=\sum_{i=1}^k u_i - \frac12 k(k+1)
\endalign
$$
(and $S_{\u}$ may be identified with a complex vector space --- or
cell --- of this dimension).

There is a similar description of the unstable manifold $U_{\u}$.
We have
$$
\lim_{t\to-\infty}
\pmatrix
e^{-a_1t} &  &  \\
  &  \ddots &  \\
  &  & e^{-a_nt}
\endpmatrix
\left[
\matrix
 &  &  \\
1 &  & \\
* & \ddots  & \\
\vdots & & 1\\
\vdots  &  & *\\
\vdots  &  &  \vdots \\
*  &  &  *
\endmatrix
\right]
=
\left[
\matrix
  &  &  \\
1 &  &  \\
  & \ddots & \\
  &  & 1\\
  &  &
\endmatrix
\right],
$$
from which we obtain a description of $U_{\u}$ of the form
$$
U_{\u}=
\pmatrix
 &  & \\
 &  & \\
1 &  & \\
0  & 1 & \\
*  & * & 0\\
0  & 0 & 1\\
*  & * & *\\
*  & * & *
\endpmatrix.
$$
This is diffeomorphic to a complex vector space whose dimension is
$$
\align
\dim_{\C} U_{\u} &=
(n-(u_1-1)-k)+(n-(u_2-1)-(k-1))+\dots+(n-(u_k-1)-1)\\
&=k(n-k) - \dim_{\C} S_{\u}\\
&=\dim_{\C} \grkcn - \dim_{\C}S_{\u}.
\endalign
$$

This identification of the stable and unstable manifolds leads to
a proof of Theorem 3.1.1 for general $k$,
as in the case $k=1$.  Note that the indices of 
the critical points are twice the complex dimensions 
of the unstable manifolds.

\no {\it Remark:} From Lemma 3.1.2, we see that the integral curves of
$-\nabla f$ \ll preserve\rr the submanifold $\grkrn$ of $\grkcn$.  In fact,
the above analysis works equally well for the restriction of $f$ to
$\grkrn$.

There is an interesting group-theoretic interpretation
of $S_{\u}$ and $U_{\u}$. The group $\glnc$ of invertible complex
$n\times n$ matrices acts naturally on $\grkcn$ (by multiplying a
column vector on the left). We have the usual exponential map
$$
\exp:\frakglnc\to\glnc,\quad
X\mapsto I+\frac{X}{1!}+\frac{X^2}{2!}+\frac{X^3}{3!}+\dots
$$
where $\frakglnc$ denotes the Lie algebra of $\glnc$, i.e. the
vector space of all $n\times n$ complex matrices. With this 
notation we see that $S_{\u}$ consists of $k$-planes of the form
$$
\exp
\pmatrix
0  & 0 & * & * & 0 & * & 0 & 0 \\
0  & 0 & * & * & 0 & * & 0 & 0 \\
0  & 0 & 0 & 0 & 0 & 0 & 0 & 0 \\
0  & 0 & 0 & 0 & 0 & 0 & 0 & 0 \\
0  & 0 & 0 & 0 & 0 & * & 0 & 0 \\
0  & 0 & 0 & 0 & 0 & 0 & 0 & 0 \\
0  & 0 & 0 & 0 & 0 & 0 & 0 & 0 \\
0  & 0 & 0 & 0 & 0 & 0 & 0 & 0 
\endpmatrix
\left[
\matrix
0  & 0 & 0 \\
0  & 0 & 0 \\
1  & 0 & 0 \\
0  & 1 & 0 \\
0  & 0 & 0 \\
0  & 0 & 1 \\
0  & 0 & 0 \\
0  & 0 & 0 
\endmatrix
\right].
$$
Thus, $S_{\u}=\{(\exp\,X)V_{\u} \st X\in \frak n_{\u} \}$, where
$\frak n_{\u}$ is a {\it nilpotent Lie subalgebra} of $\frakglnc$.
This means that $S_{\u}$ is the orbit $V_{\u}$ under the corresponding
Lie group $N_{\u}$. (Note that $\exp\,X = I + X$ here.) There is
an analogous description of $U_{\u}$.

The action of $A=\diag(e^{-a_1t},\dots,e^{-a_nt})$ (which gives the integral
curves of $-\nabla f$) is easy to express in term of the Lie algebra
$\frak n_{\u}$. For $X\in \frak n_{\u}$, we have
$$
\align
A (\exp\,X) V_{\u} &= A (\exp\,X) A^{-1} A V_{\u}\\
& = A (\exp\,X) A^{-1} V_{\u} \ \text{(as $V_{\u}$ is fixed by $A$)}\\
& = (\exp\,AX A^{-1}) V_{\u}.
\endalign
$$
The map $X\mapsto AX A^{-1}$ has the effect of multiplying 
the $(i,j)$-th entry of $X$ by $e^{-(a_i-a_j)t}$.

\subheading{3.2 Morse-Bott functions on Grassmannians}

By relaxing the condition $a_1 > \dots > a_n \ge 0$ to
$a_1 \ge \dots \ge a_n \ge 0$, we obtain a Morse-Bott function on $\grkcn$.
To investigate this, we introduce the notation
$$
\C^n = E_1 \oplus \dots \oplus E_l
$$
for the eigenspace decomposition of $D=\diag(a_1,\dots,a_n)$.
We write $b_i$ for the eigenvalue on $E_i$. Thus, the {\it distinct}
$a_i^{,}s$ are $b_1>\dots > b_l$.

\proclaim{Theorem 3.2.1} Assume that $a_1\ge\dots\ge a_n\ge 0$. Then
$f$ is a perfect Morse-Bott function. 
An element $V$ of $\grkcn$  is a critical
point of $f$ if and only if 
$V=V\cap E_1\oplus\dots\oplus V\cap E_k$,
i.e. if and only if $V$ is spanned by eigenvectors of 
$D=\diag(a_1,\dots,a_n)$. The 
NDCM containing such a $V$ is
diffeomorphic to
$Gr_{c_1}(E_1)\times\dots\times Gr_{c_l}(E_l)$
where $c_i = \dim_{\C} V\cap E_i$.
\endproclaim

It may seem unnecessary to introduce Morse-Bott functions in
a situation like this, where we already have a good supply of
Morse functions. However, we shall see later that the special
properties of Morse-Bott functions can be extremely useful.

To prove the theorem, and to identify the stable and unstable
manifolds, we use the explicit formula for the integral curves 
of $-\nabla f$ which was given earlier. 
(The derivation of this formula is clearly valid for
arbitrary real $a_1,\dots,a_n$.) We shall just sketch the main points here.

First, it is immediate from the form of the integral curves
that the critical points are as stated in the theorem. Note
that every $V_{\u}$ is certainly a critical point of $f$, but
these are not the only critical points; the others are obtained
by taking the orbits of the $V_{\u}^,$s under the product of unitary groups
$U(E_1)\times\dots\times U(E_l)$.

Let us now try to identify the stable manifold of a critical
point $V_{\u}$. Each NDCM contains at least one point of this
form. Moreover, since the negative bundles turn out to be homogeneous
vector bundles on the NDCM 
$M_{\u} \cong Gr_{k_1}(E_1)\times\dots\times Gr_{k_l}(E_l)$,
the stable manifold of $V_{\u}$ will determine the entire negative bundle.

Let $g$ be the Morse function on $\grkcn$ corresponding to
distinct real numbers $c_1>\dots > c_n\ge 0$, as in the
previous section.  The stable manifold $S^f_{\u}$ of $V_{\u}$
for $f$ is related to the stable manifold $S^g_{\u}$ of $V_{\u}$
for $g$, as we shall illustrate by considering the case
of $Gr_3(\C^8)$, with $b_1=a_1=a_2=a_3=a_4$ and
$b_2=a_5=a_6=a_7=a_8$. Recall that the stable manifold of 
$V_{\u}=V_3\oplus V_4 \oplus V_6$ (for example)
for $g$ is obtained by considering integral curves of the form:
$$
\pmatrix
e^{-c_1 t} & & & & & & &  \\
 &e^{-c_2 t} & & & & & &  \\
 & &e^{-c_3 t} & & & & &  \\
 & & &e^{-c_4 t} & & & &  \\
 & & & &e^{-c_5 t} & & &  \\
 & & & & &e^{-c_6 t} & &  \\
 & & & & & &e^{-c_7 t} &  \\
 & & & & & & &e^{-c_8 t}  
\endpmatrix
\left[
\matrix
* &* &* \\
* &* &* \\
1 &0 &0 \\
0 &1 &0 \\
 & &* \\
 & &1 \\
 & &0 \\
 & &0 
\endmatrix
\right]
$$
The NDCM $M_{\u}$ is diffeomorphic
to $Gr_2(\C^4)\times Gr_1(\C^4)$.  Let us consider the integral
curves of $-\nabla f$ through the points of $S^g_{\u}$:
$$
\pmatrix
e^{-b_1 t} & & & & & & &  \\
 &e^{-b_1 t} & & & & & &  \\
 & &e^{-b_1 t} & & & & &  \\
 & & &e^{-b_1 t} & & & &  \\
 & & & &e^{-b_2 t} & & &  \\
 & & & & &e^{-b_2 t} & &  \\
 & & & & & &e^{-b_2 t} &  \\
 & & & & & & &e^{-b_2 t}  
\endpmatrix
\left[
\matrix
* &* &* \\
* &* &* \\
1 &0 &0 \\
0 &1 &0 \\
 & &* \\
 & &1 \\
 & &0 \\
 & &0 
\endmatrix
\right]
$$
Clearly these integral curves do not necessarily approach $V_{\u}$
as $t\to\infty$; but they do if we set the components denoted below
by $\#$ equal to zero:
$$
\left[
\matrix
\# &\# &* \\
\# &\# &* \\
1 &0 &0 \\
0 &1 &0 \\
 & &\# \\
 & &1 \\
 & &0 \\
 & &0 
\endmatrix
\right]
$$
Thus, we have identified a $2$-dimensional cell in $S^f_{\u}$
which is contained in  $S^g_{\u}$.  It can be shown that this cell
is {\it precisely} $S^f_{\u}$.  In general, the procedure by which we identify
$S^f_{\u}$ as a subspace of $S^g_{\u}$ is that we delete those
coordinates which correspond to the NDCM $M_{\u}$.
It is possible to give a more systematic description of these stable
(and unstable) manifolds, as orbits of certain subgroups of $\glnc$.

\subheading{3.3 Homology groups of Grassmannians}

In section 3.1 we described a perfect Morse function $f:\grkcn\to\R$.
It follows that the homology groups of $\grkcn$ may be read off
from the Morse polynomial of $f$. For $k=1$ this is easy, but for general 
$k$ it is a little harder. In this section we shall carry out the
calculation, in two quite different ways.

Recall that $f$ has $\binom{n}{k}$ critical points, namely the
$k$-planes $V_{\u}=V_{u_1}\oplus \dots \oplus V_{u_k}$, where
$1\le {u_1}< \dots < {u_k} \le n$.  From 3.2 the (complex) dimension
of the stable manifold $S_{\u}$ is $\sum_{i=1}^k (u_i-i)$, so
the index of $V_{\u}$ as a critical point of $-f$ is
$2\sum_{i=1}^k (u_i-i)$.

If we define
$$
a_d=
\vert\{ \u \st \sum_{i=1}^k (u_i-i) = d \}\vert
$$
then the Morse polynomial of $-f$ is $M(t)=\sum_{d\ge 0} a_d t^{2d}$.
Our task, therefore,  is to calculate $a_d$.

There is a one to one correspondence between
$$
\{ (u_{1},\dots,u_{k}) \st
1 \le u_1 < \dots < u_k \le n, \ \sum_{i=1}^k (u_i-i) = d \}
$$
and
$$
\{ (p_{1},\dots,p_{k}) = (u_{1}-1,\dots,u_{k}-k) \st
0 \le p_1 \le \dots \le p_k \le n-k, \ \sum_{i=1}^k p_i = d \}.
$$
Thus, $a_d$ is equal to \ll the number of partitions of $d$ into
at most $k$ integers, where each such integer is at most $n-k$\rrr.
Unfortunately there is no simple formula for $a_d$, but there is 
a formula for the generating function
$G(t) = \sum_{d\ge 0} a_d t^{d}$ (and this is exactly what we want,
since $M(t)=G(t^2)$).

Before stating this formula, we recall the well known fact
that, if $b_d$ denotes
the number of (unrestricted) partitions of $d$, then
$$
\sum_{d\ge 0} b_d t^{d} =
\frac {1}
{(1-t)(1-t^2)(1-t^3)\dots}
$$
This means that $b_d$ is equal to the coefficient of $t^d$ in the
formal expansion of the right hand side.
The formula for $G(t)$ is:

\proclaim{Proposition 3.3.1} The generating function 
$G(t) = \sum_{d\ge 0} a_d t^{d}$ is given by
$$
G(t) = \frac 
{(1-t^{n-k+1})(1-t^{n-k+2})\dots (1-t^{n})}
{(1-t)(1-t^2)\dots (1-t^k)}
$$
\endproclaim

\no This is a purely combinatorial statement, which may be proved
directly. However,
we shall give an indirect proof, by making use of a Morse-Bott function
on $\grkcn$. In doing so, we shall find a simple inductive formula
for $G(t)$ as well.

Let $g:\grkcn\to \R$ be the Morse-Bott function corresponding to
the choice of real numbers
$$
a_1=1,\ a_2=\dots =a_n=0.
$$
Amongst the non-constant Morse-Bott functions of section 3.2, this is the
\ll crudest\rrr, in the sense that it that it has the least number
of critical manifolds. (Morse functions are at the opposite extreme;
they have the largest number of critical manifolds.) 

From Theorem 3.2.1, the critical manifolds are as follows:

\no $M_{max}=
\{ V\in\grkcn \st V = V_1\oplus W, W\sub V_{2}\oplus\dots\oplus V_{n} \}
\ \cong\  Gr_{k-1}(\C^{n-1})$

\no $M_{min}=
\{ V\in\grkcn \st V \sub V_{2}\oplus\dots\oplus V_{n} \}
\ \cong\  Gr_{k}(\C^{n-1})$.

The unstable manifold of $M_{max}$ must be the complement of $M_{min}$,
so the index of $M_{max}$ for $g$ is 
$2(\dim \grkcn - \dim Gr_{k-1}(\C^{n-1}))$, i.e. $2(n-k)$.
The index of $M_{min}$ is of course zero.
So the Morse-Bott polynomial of $g$ is
$$
MB(t) = P_{k,n-1}(t) + t^{2(n-k)}P_{k-1,n-1}(t)
$$
where $P_{i,j}(t)$ denotes the Poincar\acuteaccent e polynomial
of $Gr_{i}(\C^{j})$.

Since $g$ is a perfect Morse-Bott function, this gives an
inductive formula for $P_{k,n}(t) = MB(t)$ :

\proclaim{Proposition 3.3.2} The Poincar\acuteaccent e polynomial
$P_{k,n}(t)$ of $\grkcn$ satisfies the relation
$$
P_{k,n}(t) = P_{k,n-1}(t) + t^{2(n-k)}P_{k-1,n-1}(t).
$$
\endproclaim

\no For example:
$P_{2,4}(t) = P_{2,3}(t) + t^{4} P_{1,3}(t)=
(1+t^2+t^4) + t^4(1+t^2+t^4)
=1 + t^2 + 2t^4 + t^6 + t^8$.

We can also use this to give a proof of Proposition 3.3.1, which
is equivalent to the following slightly more symmetrical statement:

\proclaim{Proposition 3.3.3} The Poincar\acuteaccent e polynomial
$P_{k,n}(t)$ of $\grkcn$ is given by
$$
P_{k,n}(t) = \frac
{\prod_{i=1}^n(1-t^{2i})}
{\prod_{i=1}^k(1-t^{2i})  \prod_{i=1}^{n-k}(1-t^{2i})}\ .
$$
\endproclaim

\demo{Proof}  Denote the right hand side by $B_{k,n}(t)$.
It is easy to verify that this satisfies the same
recurrence relation as $P_{k,n}(t)$, i.e.
$B_{k,n}(t) = B_{k,n-1}(t) + t^{2(n-k)}B_{k-1,n-1}(t)$.
Since the recurrence relation determines $B_{k,n}(t)$
or $P_{k,n}(t)$ inductively,  and these agree when $k=1$,
they must be equal.
\qed\enddemo

\subheading{3.4 Schubert varieties}

In section 3.1 we described explicitly the stable manifold $S_{\u}$
of a critical point $V_{\u}=V_{u_1}\oplus\dots\oplus V_{u_k}$
of a Morse function $f:\grkcn\to \R$.  It has the form
$$
S_{\u}=
\left\{
\left.
V=
\left[
\matrix
* & * & *\\
* & * & *\\
1 & 0 & 0\\
  & 1 & 0\\
  &  & *\\
  &  & 1\\
  &  & 0\\
  &  & 0
\endmatrix
\right]
\in\grkcn
\ \right|\ 
\ast\in\C
\right\}.
$$
From 3.1 it is clear that such $V$ are characterized {\it geometrically}
by the following conditions:
$$
\dim V\cap \C^i =
\cases
0 \quad\text{if}\quad 1\le i \le u_1-1\\
1 \quad\text{if}\quad  u_1\le i \le u_2-1\\
\quad\quad\dots\\
k\quad\text{if}\quad  u_k\le i \le n
\endcases
$$
In other words, we have
$$
S_{\u}=
\{ V\in\grkcn \st
\dim V\cap \C^i = \dim V_{\u}\cap \C^i \ \text{for all}\ i \}.
$$
The condition on $\dim V\cap \C^i$ is called a \ll Schubert condition\rrr.
It can be specified either by listing 
$v_i=\dim V\cap \C^i$,
i.e.
$$
v_1=\dots=v_{u_1-1}=0,\ 
v_{u_1}=\dots=v_{u_2-1}=1, \ \  \dots, \ \  
v_{u_k}=\dots=v_n=k
$$
or, more economically, by listing those $i$ such that
$\dim V\cap \C^i = \dim V\cap \C^{i-1} + 1$, i.e.
$$
u_1,u_2,\dots,u_k.
$$
The $k$-tuple $\u=(u_1,\dots,u_k)$ is referred to as a
\ll Schubert symbol\rrr, and the set $S_{\u}$ is called the
\ll Schubert cell\rr associated to $\u$.

Similarly, the unstable manifold $U_{\u}$ is characterized
by these conditions:
$$
\dim V\cap (\C^{n-i})^\perp =
\cases
0 \quad\text{if}\quad 1\le i \le n- u_k\\
1 \quad\text{if}\quad  n- u_k +1\le i \le n- u_{k-1}\\
\quad\quad\dots\\
k\quad\text{if}\quad  n- u_1 +1\le i \le n
\endcases
$$

\proclaim{Example 3.4.1}
\rm 
For the Morse function $f:Gr_2(\C^4)\to \R$ (of section 3.1), there
are $6$ Schubert cells. We give the matrix representations below, followed
by the sequence $\dim V\cap \C,\dots,\dim V\cap \C^4$, the Schubert
symbol, and the dimension of the cell.
$$
\pmatrix
* &* \\
* &* \\
1 &0 \\
0 &1 
\endpmatrix
\quad\quad
0,0,1,2
\quad\quad
(3,4)
\quad\quad
\dim_{\C}=4
$$
{}
$$
\pmatrix
* &* \\
1 &0 \\
0 &* \\
0 &1 
\endpmatrix
\quad\quad
0,1,1,2
\quad\quad
(2,4)
\quad\quad
\dim_{\C}=3
$$
{}
$$
\pmatrix
* &* \\
1 &0 \\
0 &1 \\
0 &0 
\endpmatrix
\quad\quad
0,1,2,2
\quad\quad
(2,3)
\quad\quad
\dim_{\C}=2
$$
{}
$$
\pmatrix
1 &0 \\
0 &* \\
0 &* \\
0 &1 
\endpmatrix
\quad\quad
1,1,1,2
\quad\quad
(1,4)
\quad\quad
\dim_{\C}=2
$$
{}
$$
\pmatrix
1 &0 \\
0 &* \\
0 &1 \\
0 &0 
\endpmatrix
\quad\quad
1,1,2,2
\quad\quad
(1,3)
\quad\quad
\dim_{\C}=1
$$
{}
$$
\pmatrix
1 &0 \\
0 &1 \\
0 &0 \\
0 &0 
\endpmatrix
\quad\quad
1,2,2,2
\quad\quad
(1,2)
\quad\quad
\dim_{\C}=0
$$
It follows that the Poincar\acuteaccent e polynomial of $Gr_2(\C^4)$
is $1+t^2+2t^4+t^6+t^8$.
\qed\endproclaim

\proclaim{Definition 3.4.2} The {\it Schubert variety} $X_{\u}$
associated to the Schubert symbol $\u$ is the closure of  $S_{\u}$
(with respect to the usual topology of $\grkcn$), i.e.
$$
X_{\u}= \overline{S}_{\u} =
\{ V\in\grkcn \st
\dim V\cap \C^i \ge v_i \ \text{\rm for all}\ i \}.
$$
\endproclaim

\no It is easy to show that $X_{\u}$ is an algebraic subvariety of
$\grkcn$ (see section 3.5).  This subvariety may have singularities.  For example,
in the case of $Gr_2(\C^4)$, we have
$$
X_{(2,4)}=\{ V\in Gr_2(\C^4) \st
\dim V\cap \C^2 \ge 1 \}
$$
(the $v_i^,s$ are given by
$(v_1,v_2,v_3,v_4)=(0,1,1,2)$, but
the conditions $\dim V\cap \C \ge 0$, $\dim V\cap \C^3 \ge 1$,
$\dim V\cap \C^4 \ge 2$ are 
automatically\footnote{If $V$, $W$ are 
linear subspaces of $\C^n$ of dimensions
$k,l$ respectively, then we have
$W/W\cap V \cong W+V/V$, and hence
$\dim W\cap V + \dim (W+V) = k+l$.}
satisfied).  The point $V=\C^2$ is a singular point of $X_{(2,4)}$,
but $X_{(2,4)} - \{ \C^2 \}$ is smooth, having the structure of
a fibre bundle over $\C P^1$ with fibre $\C P^2 - \{  \text{point} \}$.

Observe that
$$
X_{(2,4)} = S_{(2,4)} \cup S_{(2,3)} \cup S_{(1,4)} \cup S_{(1,3)} \cup
S_{(1,2)}.
$$
It is clear from the definition that, in general, $X_{\u}$ is
a disjoint union of Schubert cells.  This gives rise to a partial
order on the set of Schubert symbols:  we define $\u_1 \le \u_2$
if and only if $\overline{S}_{\u_1} \supseteq S_{\u_2}$.

In the case of $ Gr_2(\C^4)$, this partial order is represented
in the following diagrams.  The second diagram indicates the
conditions which define the Schubert varieties.
$$
\quad
\matrix
& & & & \\
& & & & \\
  &  &(1,2)  &  &  \\
  &  &  &  &  \\
  &  &(1,3)  &  &  \\
  &  &  &  &  \\
(2,3)  &  &  &  &(1,4)  \\
  &  &  &  &  \\
  &  &(2,4)  &  &  \\
  &  &  &  &  \\
  &  &(3,4)  &  & \\
& & & &\\
& & & & 
\endmatrix
\quad\quad\quad\quad
\matrix
& & & & \\
& & & & \\
  &  &{V=\C^2}  &  &  \\
  &  &  &  &  \\
  &  &{\C\sub V\sub \C^3}  &  &  \\
  &  &  &  &  \\
{V\sub \C^3}  &  &  &  &{\C \sub V}  \\
  &  &  &  &  \\
  &  &{\dim V\cap\C^2\ge 1}  &  &  \\
  &  &  &  &  \\
  &  &\text{no condition}  &  & \\
& & & &\\
& & & & 
\endmatrix
$$
Although the partial order is a simple consequence of the
definition of Schubert variety, it provides nontrivial information
on the behaviour of the flow lines of $-\nabla f$.  Namely, the
condition $\u_1\le \u_2$ is equivalent to the existence of a
flow line $\ga(t)$ from $\u_2$ to $\u_1$, i.e.
such that $\lim_{t\to -\infty}\ga(t)=V_{\u_2}$,
$\lim_{t\to \infty}\ga(t)=V_{\u_1}$. (This is not immediately
obvious from the condition $\overline{S}_{\u_1} \supseteq S_{\u_2}$,
but it does follow from the
geometrical description of $S_{\u_1}\cap U_{\u_2}$.)

We shall return to Schubert varieties in section 3.6, when we discuss
the cohomology ring of a Grassmannian, so we conclude this section with 
some further remarks.

First, since the flow lines of $-\nabla f$ \ll preserve\rr the
Schubert variety $X_{\u}$, we deduce that $X_{\u}$ inherits a decomposition
into (possibly singular) \ll stable manifolds\rr 
(or \ll unstable manifolds\rrr).

Second, although the Morse function $f$ depends on a choice of 
real numbers $a_1>\dots>a_n$, the Schubert varieties are independent
of this choice.  In fact they depend solely on the standard \ll flag\rr
$$
\C \sub \C^2 \sub \dots \sub \C^n,
$$
or on the standard ordered basis 
$e_1,\dots,e_n$ of $\C^n$.  If $a_1,\dots,a_n$ are arbitrary distinct real
numbers, we obtain a similar Morse function, possibly corresponding
to a re-ordering of $e_1,\dots,e_n$.  More generally still, any
choice of orthonormal basis, or equivalently any flag, 
corresponds to a similar
Morse function.  The formula for such a function is obtained by replacing
the diagonal matrix $A$ by $UAU^{-1}$, where $U$ is a unitary matrix.

There is a geometrical description of the stable manifolds of the
Morse-Bott functions considered in section 3.2, i.e. where the
real numbers $a_1\ge\dots\ge a_n$ are not necessarily distinct.  We state
the result without proof (as the easiest proof depends on the more
general theory of \S 4).  Recall that the critical manifolds are
denoted $M_{\u}$, and that each such NDCM contains at least one point
of the form $V_{\u}$. The stable manifold $S_{\u}$ of $M_{\u}$, i.e. the
union of the stable manifolds of all points of $M_{\u}$, is
then given by:
$$
S_{\u}=
\{ V\in\grkcn \st
\dim V\cap (E_1\oplus\dots\oplus E_i) = 
\dim V_{\u}\cap (E_1\oplus\dots\oplus E_i)
\ \text{for all}\ i \},
$$
where $E_1,\dots,E_l$ are the eigenspaces of $\diag(a_1,\dots,a_n)$
as in section 3.2.
Thus, in this case, the Schubert \ll cells\rr (or rather, Schubert
cell-bundles) depend only on the \ll partial flag\rr
$$
E_1 \sub E_1\oplus E_2 \sub \dots \sub E_1\oplus\dots\oplus E_l = \C^n.
$$
Conversely, as in the case of the Morse functions discussed earlier,
{\it any} partial flag determines a Morse-Bott function on $\grkcn$.

The Schubert cell-bundles (or their NDCMs) are parametrized by 
$l$-tuples $(c_1,\dots,c_l)$  of non-negative integers with
$c_j\le \dim E_j$ and $c_1+\dots+c_l=k$; namely
$c_j=\dim V_{\u}\cap E_j$.  We consider  $(c_1,\dots,c_l)$
to be a \ll generalized Schubert symbol\rrr. The bundle projection
map $S_{\u}\to M_{\u}$ is given explicitly by
$$
V\mapsto (V_{(1)},\dots,V_{(l)})\in
Gr_{c_1}(E_1)\times\dots\times Gr_{c_l}(E_l)
$$
where
$$
V_{(i)} = V\cap E_1\oplus\dots\oplus E_i + E_1\oplus\dots\oplus E_{i-1}\ /
\ E_1\oplus\dots\oplus E_{i-1}.
$$

The integers $w_i=c_1+\dots+c_i$, $i=1,\dots,l$,
are analogous to the integers $v_i$ in the case of a Morse
function.  In terms of these integers, the Schubert cell-bundle 
$S_{\u}$ is
$$
S_{\u}=
\{ V\in\grkcn \st
\dim V\cap (E_1\oplus\dots\oplus E_i) = w_i
\ \text{for all}\ i \},
$$
By taking the closure of a Schubert cell-bundle, we obtain a generalized
Schubert variety, namely
$$
X_{\u}= \overline{S}_{\u} =
\{ V\in\grkcn \st
\dim V\cap E_1\oplus\dots\oplus E_i \ge w_i
\ \text{for all}\ i \}.
$$
One of the advantages of having these explicit descriptions
of Schubert cell-bundles is that we can compute easily the indices
of the critical points; we shall give some examples below.

\proclaim{Example 3.4.3}
\rm
Let us choose the partial flag $\C\sub \C^n$. This corresponds to a
Morse-Bott function on $\grkcn$ 
with $a_1>a_2=\dots=a_n$. We
have already considered such a function in section 3.3; there are
two critical manifolds.  The stable manifold of the maximum NDCM is
just that NDCM, and the stable manifold of the minimum NDCM is
the complement of the maximum NDCM.  

The eigenspace decomposition of $\C^n$ is given by $E_1=\C$,
$E_2=\C^\perp$, and the generalized Schubert symbols are
$(c_0,c_1)=(1,k-1)$ and $(c_0,c_1)=(0,k)$.
\qed\endproclaim

\proclaim{Example 3.4.4}
\rm
Consider the Morse-Bott function on $\C P^6 = Gr_1(\C^7)$ given
by $a_1=a_2=2$, $a_3=a_4=a_5=1$, $a_6=a_7=0$.  In this case the
eigenspace decomposition of $\C^7$ is given by 
$$
\C^7=E_1\oplus E_2\oplus E_3, \quad
E_1=V_1\oplus V_2,
E_2=V_3\oplus V_4\oplus V_5, E_3=V_6\oplus V_7.
$$  
We list below the
generalized Schubert symbols, followed by the NDCM, and then the
Schubert cell-bundle.

\no $(c_0,c_1,c_2)=(1,0,0)$,
$M_{\u}=\P(E_1)$,
$S_{\u}=\P(E_1)$

\no $(c_0,c_1,c_2)=(0,1,0)$,
$M_{\u}=\P(E_2)$,
$S_{\u}=\P(E_1\oplus E_2)-\P(E_1)$

\no $(c_0,c_1,c_2)=(0,0,1)$,
$M_{\u}=\P(E_3)$,
$S_{\u}=\P(E_1\oplus E_2\oplus E_3)-\P(E_1\oplus E_2)$

In particular, the Morse indices (for the function $-f$) are
$0, 4, 10$ respectively, and the Morse-Bott polynomial is
$t^0(1+t^2) + t^4(1+t^2+t^4) + t^{10}(1+t^2)$. This is equal
to the Poincar\acuteaccent e polynomial of $\C P^6$, as
expected.
\qed\endproclaim

\subheading{3.5 Morse theory of the Pl\"ucker embedding}

There is a well known embedding
$$
\grkcn \to \C P^{N},\quad N = \binom{n}{k} - 1
$$
called the  Pl\"ucker embedding.  It is defined by
$$
V\mapsto \wedge^k V \sub \wedge^k\C^n \cong \C^{N+1}.
$$
If $e_1,\dots,e_n$ are the standard basis vectors of $\C^n$,
then the vectors $e_{\u}=e_{u_1}\wedge\dots\wedge e_{u_k}$
with $1\le u_1 < \dots < u_k \le n$
constitute a basis of $\wedge^k\C^n$.

Let $f:\grkcn\to\R$ be the Morse-Bott function defined by
certain real numbers $a_1,\dots,a_n$, as in section 3.2. Then
by Lemma 3.1.2
the one parameter diffeomorphism group of the vector field
$-\nabla f$ is induced by the action
$$
t\cdot \sum \la_i e_i =  \sum \la_i e^{-ta_i} e_i
$$
of $\R$ on $\C^n$. 

This action is the restriction of the action 
$$
t\cdot \sum \la_{\u} e_{\u} =  
\sum \la_{\u} e^{-t(a_{u_1}+\dots+a_{u_k})} e_{\u}
$$
of $\R$ on $\wedge^k\C^n$.  But this action is the one parameter 
diffeomorphism group of the vector field $-\nabla F$, where 
$F:\C P^{N}\to\R$ is the Morse-Bott function defined by the
$\binom{n}{k}$
real numbers $a_{u_1}+\dots+a_{u_k}$. 
We conclude that the Morse-Bott theory of $f$ on $\grkcn$ is just
the \ll restriction\rr of the (much simpler!) Morse-Bott theory
of $F$ on $\C P^{N}$:

\proclaim{Proposition 3.5.1}
Let $M_{\u},S_{\u},X_{\u}$ denote the NDCM's, Schubert cell-bundles,
and generalized Schubert varieties for the Morse-Bott function
$f:\grkcn\to\R$.  Let $M_{\u}^{F},S_{\u}^{F},X_{\u}^{F}$ denote
the corresponding objects for the Morse-Bott function
$F:\C P^{N}\to\R$. Then we have
$M_{\u} = \grkcn \cap M_{\u}^{F}$,  $S_{\u} = \grkcn \cap S_{\u}^{F}$,
and $X_{\u} \sub \grkcn \cap X_{\u}^{F}$.
\qed\endproclaim

\no Observe that it is possible to choose the real numbers 
$a_1,\dots,a_n$ so that both $f$ and $F$ are Morse functions.
But it may happen that $f$ is a Morse function even when $F$
is not.

The spaces $M_{\u}^{F},S_{\u}^{F},X_{\u}^{F}$ are determined 
by an eigenspace decomposition 
$$
\C^{N+1}= \hat E_1 \oplus \dots
\oplus \hat E_{\hat k}
$$ 
in the usual way.  Using this notation,
the NDCM's of $F$ are the linear subspaces $\P(\hat E_i)$, and the
stable manifold of $\P(\hat E_i)$ is given explicitly as a bundle over 
$\P(\hat E_i)$ by
$$
\P(\hat E_0\oplus\dots\oplus \hat E_i) - \P(\hat E_0\oplus\dots\oplus \hat E_{i-1}).
$$
The projection map to $\P(\hat E_i)$ sends a line $L$ in 
$\hat E_0\oplus\dots\oplus \hat E_i$ to the line 
$L + \hat E_0\oplus\dots\oplus \hat E_{i-1}$ in 
$\hat E_0\oplus\dots\oplus \hat E_i\ /\ \hat E_0\oplus\dots\oplus \hat E_{i-1}
\cong \hat E_i$.

The associated Schubert variety, i.e. the closure of this
stable manifold, is just
$$
\P(\hat E_0\oplus\dots\oplus \hat E_i),
$$
which is a {\it linear subspace} of $\C P^{N}$.  Although $X_{\u}$
is not necessarily {\it equal} to the intersection of this space
with $\grkcn$, it is in fact true that  $X_{\u}$
is given by the intersection of {\it some} linear subspace
with $\grkcn$.  This follows from the fact that
$$
S_{\u}=\grkcn \quad \cap \quad
\P(\hat E_0\oplus\dots\oplus \hat E_i) - \P(\hat E_0\oplus\dots\oplus \hat E_{i-1})
$$
(as some of the linear equations defining 
$\P(\hat E_0\oplus\dots\oplus \hat E_i)$ in $\C^{n+1}$
may become dependent in the presence of the equations
defining $\grkcn$).

The fact that the Pl\"ucker embedding is compatible with the
natural Morse-Bott functions on $\grkcn$ and $\C P^{N}$ may be
explained group-theoretically.  The key point is that the 
Pl\"ucker embedding is induced by an irreducible representation
$U_n\to U_{N+1}$.  However, it seems technically easier to work
with the explicit formulae for the flow lines, 
as we have done in this section.

\subheading{3.6 Cohomology of the Grassmannian, and the Schubert calculus}

In this section we consider only Morse functions on $\grkcn$.

From a Schubert symbol $\u=(u_1,\dots,u_k)$ we obtain a
cell $S_{\u}$ in $\grkcn$ of (real) dimension $2\sum_{i=1}^k (u_i - i)$.
By the Morse inequalities for the coefficient group $\Z$ (or by
standard theory of cellular homology) it follows that
$H_\ast(\grkcn;\Z)$ is a free abelian group with one generator
for each Schubert cell. Let $x_{\u}$ be the homology class
represented by\footnote{The precise meaning of 
this is explained in Appendix B
of \cite{Fu1}.}
$X_{\u}$.  Both $x_{\u}$ and $X_{\u}$ are referred
to as \ll Schubert cycles\rrr.

By Poincar\acuteaccent e Duality we  obtain a dual cohomology class  $z_{\u}$
of dimension $2k(n-k) - 2\sum_{i=1}^k (u_i - i)$. Thus,
$\dim z_{\u} = \codim x_{\u}$.

An example of particular interest is the generator of $H^2(\grkcn;\Z)\cong\Z$;
this corresponds to the unique codimension one Schubert cycle, 
i.e. to the Schubert symbol $(n-k,n-k+2,\dots,n-1,n)$. The Schubert
conditions here are
$$
\dim V\cap\C = \dots = \dim V\cap\C^{n-k-1}=0,\ 
\dim V\cap\C^{n-k}=1,\  \dots, \ 
\dim V\cap\C^n=k,
$$
and the Schubert variety is characterized by the single condition
$\dim V\cap\C^{n-k}\ge1$. In terms of the Pl\"ucker embedding
(see section 3.5), we have
$$
X_{(n-k,n-k+2,\dots,n-1,n)}=\grkcn \cap H
$$
where $\P(H)$ is the Schubert variety for the critical
point $V_{n-k}\wedge V_{n-k+2}\wedge \dots \wedge V_{n-1} \wedge V_{n}$
in $\C P^N$. Now, since $a_1>\dots > a_n$, we have 
$a_{n}+\dots + a_{n-k+1} > a_{n} + \dots + a_{n-k+2} + a_{n-k} > \dots$, 
so $H$ is the {\it hyperplane} in $\C^{N+1}$ orthogonal to
$V_{n-k+1}\wedge \dots \wedge V_{n-1} \wedge V_{n}$.
It is well known that the cohomology class dual to $\P(H)$ is
a generator of $H^2(\C P^N;\Z)$, so we deduce that the induced
homomorphism $H^2(\C P^N;\Z) \to H^2(\grkcn;\Z)$ is an isomorphism.

The muliplicative behaviour of $H^\ast(\grkcn;\Z)$ is equivalent
to the behaviour of the intersections of generic representatives
of homology classes (we shall make a more precise statement
shortly). As a first step towards describing this, 
we shall need a slight generalization of Schubert varieties.

In section 3.4 we pointed out that the definition of the Schubert
varieties $X_{\u}$ depends only on the choice of the standard flag
$\C \sub \C^2 \sub \dots \sub \C^n$. If we choose a new flag
$F_1 \sub F_2 \sub \dots \sub F_n = \C^n$, denoted by $F$,
then we obtain new
objects $S_{\u}^F,X_{\u}^F,x_{\u}^F,z_{\u}^F$ defined in exactly the
same way as $S_{\u},X_{\u},x_{\u},z_{\u}$, but using the new
flag instead of the standard flag. Since any two flags are related
by an element of the unitary group $U_n$, however,
we have $x_{\u}^F=x_{\u}$ and  $z_{\u}^F=z_{\u}$. So we may regard
the Schubert cycles $X_{\u}^F$ as a family of representatives
for the homology class $x_{\u}$, parametrized by the space
of all flags.

For example, consider the \ll opposite\rr flag
$$
(\C^{n-1})^\perp \sub (\C^{n-2})^\perp \sub \dots \sub \C^\perp \sub \C^n;
$$
let us denote the Schubert varieties with respect to this flag
by $X_{\u}^c$. It is easy to check that
$$
X_{{\u}}^c \ = \overline U_{\u^c}
$$
where $\u^c = (n-u_k +1, n-u_{k-1}+1, \dots, n-u_1+1)$.
Thus, both  $\overline S_{\u}$ and  $\overline U_{\u^c}$ 
are representatives of the same homology class  $x_{\u}$.

We now state a special case of an important general principle
(see Appendix B of \cite{Fu1}):

\proclaim{Theorem 3.6.1} If $x_{\u}$ and $x_{\v}$ are Schubert
cycles with $\dim z_{\u} + \dim z_{\v} = \dim \grkcn$, then
the product $z_{\u}z_{\v} \in H^{\dim \grkcn}(\grkcn;\Z) \cong \Z$ 
of the corresponding cohomology classes is equal to the
intersection number of $X_{\u}$ and $X_{\v}$.  
\qed\endproclaim

\no This intersection number is equal to the number of points 
(counted with multiplicities) in
$X_{\u}^{F_1} \cap X_{\v}^{F_2}$
whenever $X_{\u}^{F_1} \cap X_{\v}^{F_2}$ is a finite set.

\proclaim{Example 3.6.2}
\rm
For any $\u$, we have $z_{\u}z_{\u^c}=1$.  This is because
the dual homology classes are represented by  
$\overline S_{\u}$ and $\overline U_{\u}$ respectively, 
and these intersect
at precisely one point, namely the critical point $V_{\u}$.
(Of course, these homology classes may also be represented
by $\overline S_{\u}$ and  $\overline S_{\u^c}$. But this
is of no interest to us as these
cycles intersect at infinitely many points.)
To see that the multiplicity of the intersection point is $1$, one can
use the Pl\"ucker embedding --- it follows from our discussion
in section 3.5 that the multiplicity of {\it any} isolated point
of intersection of two Schubert varieties is precisely $1$,
since we are just taking the intersection of linear subspaces
in $\C P^N$.
\qed\endproclaim

More generally, we have:

\proclaim{Proposition 3.6.3} If $x_{\u}$ and $x_{\v}$ are Schubert
cycles with $\dim z_{\u} + \dim z_{\v} = \dim \grkcn$, then
$$
z_{\u}z_{\v} =
\cases
1  &\text{if}\ \u=\v^c\\
0 &\text{otherwise}
\endcases
$$
\endproclaim

\demo{Proof} We have just seen that $z_{\u}z_{\u^c}=1$. To show that
$X_{\u}\cap X_{\v}=\es$ if $\v\ne \u^c$, one may use the geometrical
characterization of $X_{\u}$ and $X_{\v}$ --- we omit the details.
\qed\enddemo

This proposition is a manifestation of Poincar\acuteaccent e Duality,
and it allows us to determine the products
$z_{\u}z_{\v}$ for {\it arbitrary} $\u,\v$.  For we may express
$z_{\u}z_{\v}$ in terms of the additive Schubert basis as
$$
z_{\u}z_{\v} = a_1 z_{\u_{(1)}} + \dots + a_r z_{\u_{(r)}}
$$
for some integers $a_1,\dots,a_r$, where 
$\dim z_{\u} + \dim z_{\v} = \dim z_{\u_{(i)}}$ for each $i$.
Then we obtain $a_i$ by multiplying the above expression by 
$z_{\u_{(i)}^c}$:
$$
z_{\u}z_{\v}z_{\u_{(i)}^c} = a_i
$$
(all other products vanish, by Proposition 3.6.3).  Thus, we
have to calculate all triple products $z_{\u}z_{\v}z_{\w}$
such that $\dim z_{\u} + \dim z_{\v}  + \dim z_{\w} = 
\dim \grkcn$. Theorem 3.6.1 generalizes to this situation, so we
have to calculate the corresponding triple intersections of
Schubert varieties.

\proclaim{Example 3.6.4}
\rm
We shall carry out the calculation of some triple products for
$Gr_2(\C^4)$, and hence determine the multiplicative structure of
$H^\ast(Gr_2(\C^4);\Z)$. First we list the additive generators:
$$
\align
z_{(3,4)} &\in H^0\\
z_{(2,4)} &\in H^2\\
z_{(1,4)}, z_{(2,3)} &\in H^4\\
z_{(1,3)} &\in H^6\\
1=z_{(1,2)} &\in H^8
\endalign
$$
Proposition 3.6.3 gives the following products:
$$
z_{(2,4)}z_{(1,3)}=z_{(1,4)}z_{(1,4)}=z_{(2,3)}z_{(2,3)}=1,
\ z_{(1,4)}z_{(2,3)}=0.
$$
Let us now try to compute $z_{(1,4)}z_{(2,4)}z_{(2,4)}$.
We must find suitably generic representing cycles $X_{\u}^F$ for these
classes, by choosing suitably generic flags $F$.

For $z_{(2,4)}$ we need two modifications of the standard
representative
$$
X_{(2,4)}= \{ V\in Gr_2(\C^4) \st \dim V \cap \C^2 \ge 1 \}.
$$
We shall choose the flags
$$
\align
F^\pr: \quad
& V_1 \sub V_1\oplus V_4 \sub V_1\oplus V_2\oplus V_4 \sub \C^4\\
F^{\pr\pr}: \quad
& V_2 \sub V_2\oplus V_4 \sub V_1\oplus V_2\oplus V_4 \sub \C^4
\endalign
$$
The corresponding cycles are:
$$
\align
X_{(2,4)}^\pr&= \{ V\in Gr_2(\C^4) \st \dim V \cap V_1\oplus V_4 \ge 1 \}\\
X_{(2,4)}^{\pr\pr}&= \{ V\in Gr_2(\C^4) \st \dim V \cap V_2\oplus V_4 \ge 1 \}
\endalign
$$
For $z_{(1,4)}$ we shall choose the flag
$$
F^{\pr\pr\pr}: \quad
V_3 \sub V_2\oplus V_3 \sub V_1\oplus V_2\oplus V_3 \sub \C^4
$$
i.e. we choose
$$
X^{\pr\pr\pr}_{(1,4)}= \{ V\in Gr_2(\C^4) \st V_3 \sub V \}.
$$
It may now be verified that
$$
X_{(2,4)}^{\pr}\cap X_{(2,4)}^{\pr\pr}\cap X^{\pr\pr\pr}_{(1,4)}
= \{ V_3\oplus V_4 \}
$$
i.e. a single point.  As in Example 3.6.2, we can see that the
multiplicity of this point is $1$. So we conclude that
$z_{(1,4)}z_{(2,4)}z_{(2,4)}=1$.
Exactly the same argument gives $z_{(2,3)}z_{(2,4)}z_{(2,4)}=1$.

The remaining (double) products in $H^\ast(Gr_2(\C^4);\Z)$ are
$$
\align
z_{(2,4)}z_{(2,4)} &= a z_{(2,3)} + b z_{(1,4)}\\
z_{(1,4)}z_{(2,4)} &= c z_{(1,3)}\\
z_{(2,3)}z_{(2,4)} &= d z_{(1,3)}.
\endalign
$$
Using the two triple products which we have just calculated, we
find that $a=b=c=d=1$.
\qed\endproclaim

The same method works for $H^\ast(\grkcn;\Z)$, although this
situation is of course more complicated.  
There are famous general formulae
for the multiplicative structure, which constitute the \ll Schubert
calculus\rrr. An elementary approach to these formulae and their
traditional applications can be found in \cite{Kl-La}; other
versions can be found in \cite{Gr-Ha}, 
\cite{Hl}, and \cite{Fu1}.

From the theory of Chern classes, there is a well known \ll closed
formula\rr for the ring structure of $H^\ast(\grkcn;\Z)$, namely
$$
\frac{\Z[c_1,\dots,c_{n-k},d_1,\dots,d_k]}
{(1+c_1+\dots+c_{n-k})(1+d_1+\dots+d_{k})=1}
$$
where $c_i,d_i\in H^{2i}(\grkcn;\Z)$.
(This is explained in detail in \S 23 of \cite{Bo-Tu};
another good reference is \cite{Mi-St}.)  It may be
checked that this agrees with our description in the 
case of $Gr_2(\C^4)$.  The dual of cohomology class $d_i$ is
represented by the Schubert variety $X_{\u}$ with
$\u=(n-k,n-k+1,\dots,n-k+i-1,n-k+i+1,\dots,n-1,n)$;
these are called \ll special Schubert varieties\rrr.

\subheading{3.7 Next steps}

We have now covered the \ll classical\rr aspects of Morse
theory, and in \S 4 we shall turn to more recent developments.
As motivation for this, we mention here a couple of points which 
arise from our study of Grassmannians. 

An immediate question is: when can the cohomology ring be determined
directly from a Morse function? It is a well known limitation of
classical Morse theory that the Morse inequalities give information
only about the {\it additive} structure of the cohomology ring.
However, we have seen that the cohomology ring of $\grkcn$ can be found
from explicit knowledge of the stable manifolds of a suitable Morse
function.  Was this just a special trick, or is there perhaps
a more general theory which works for Morse functions on arbitrary
compact manifolds?

A slightly more subtle (but related) question concerns the
possible configurations of flow lines of a Morse function. We saw in \S 2
that the possible configurations of critical points of a Morse
function $f:M\to\R$  are limited by
the topology of $M$. For example, it is not possible to have
a Morse function on $S^1\times S^1$ with precisely three
critical points whose indices are $0, 1$ and $2$.  In the same
way, it is not possible to have arbitrary configurations of
flow lines connecting the critical points.  The question arises
as to how these configurations of flow lines are restricted.

These questions may be answered by generalizing Morse theory
in a new way: instead of considering a single Morse function, we
consider several. Indeed, when we computed triple products of cohomology
classes in section 3.6, we were in fact making use of three
\ll independent\rr Morse functions on $\grkcn$. The general theory 
underlying this calculation is what we shall study next.

$${}$$
\head
\S 4.  Morse theory in the 1990's
\endhead

\no{\it Brief summary:}  We describe several recent 
applications of Morse theory, in which the gradient flow lines
play a fundamental role.  Although the level of discussion will be
somewhat more advanced in this section, the case of a complex
Grassmannian should be kept in mind as a typical example.
We begin by discussing Morse functions which arise from
torus actions on K\"ahler manifolds; these functions,
which include the functions on Grassmannians in \S 3, have
the crucial property that their gradient flow lines are
explicitly identifiable.  Then we describe the \ll new\rr approach to
Morse theory, due to Witten. After that we present the
\ll field-theoretic\rr Morse theory of Cohen-Jones-Segal,
Betz-Cohen, and Fukaya.

\subheading{4.1 Morse functions generated by torus actions}

In \S 1 and \S 2 we gave a review of the \ll classical\rr
Morse theory, and then in \S 3 we illustrated this in detail for a
particular example.  We shall now focus on
some contemporary aspects, which show that the power of Morse
theory goes far beyond the computation of homology groups.
Our starting point is an important family of examples
which includes the Morse and Morse-Bott functions of \S 3.

Let $T\cong S^1\times\dots\times S^1$ be a torus, and let
$M$ be a simply-connected connected compact K\"ahler manifold.
Assume that the group $T$ acts on $M$, and that this action preserves
the complex structure $J$ and the K\"ahler $2$-form $\om$ of $M$.
It follows that the action also preserves the  K\"ahler
metric $\lan\ ,\ \ran$, as the latter is given by
$\om(A,B)=\lan A,JB \ran$.

Let $\t$ denote the Lie algebra of $T$. Let $X\in \t$
be any generator of the torus; this means that $T$ is the closure of
its (not necessarily closed) subgroup $\exp\,\R X$. 
The formula
$$
X^\ast_m = \ddt \exp\, tX \cdot m \vert_{t=0}
$$
defines a vector field $X^\ast$ on $M$.  

Since $M$ is simply connected, 
and $\om(X^\ast, \ )$ is a closed $1$-form,
there is a function $f^X:M\to \R$
such that $df^X=\om(X^\ast, \ )$. We shall see
that $f^X$ is a perfect Morse-Bott function, with particularly
nice properties.  Observe that we have
$$
-\nabla f^X = JX^\ast,
$$
from the formula $\lan \nabla f^X, A\ran = df^X(A) = \om(X^\ast,A)
= - \lan JX^\ast,A \ran$.

\proclaim{Theorem 4.1.1 (\cite{Fr1})} For
any generator $X$ of $\t$, the function $f^X$ is a perfect
Morse-Bott function on $M$. The critical points of $f^X$ are the fixed
points of the $T$-action, and they form a finite number
of totally geodesic K\"ahler submanifolds of $M$.
\endproclaim

\demo{Sketch of the proof} The fact that the critical points of
$f^X$ are the fixed points of $T$ follows from the formula 
$-\nabla f^X = JX^\ast$.

Let $m$ be a fixed point of $T$. Then
there is an induced action of the Lie group $T$, and hence
also of the Lie algebra $\t$, on the vector space $T_mM$.
This means that we have a Lie group homomorphism
$\Th_m:T\to Gl(T_m M)$ and a Lie algebra homomorphism
$\th_m:\t \to \End(T_m M)$.

Since $T$ acts by isometries (i.e. the
action of $T$ preserves the metric), $\Th_m(t)$ is an
orthogonal transformation, and $\th_m(X)$ a skew-symmetric
transformation, for each $t\in T$, $X\in \t$.
As $T$ is abelian, we may put these transformations
simultaneously into canonical form.  This means that
there exists a decomposition
$$
T_m M = V_0 \oplus V_1 \oplus \dots \oplus V_r
$$
such that $\th \vert_{V_0} = 0$ and
$$
\th\vert_{V_i} =
\pmatrix
0 & w_i \\
-w_i & 0
\endpmatrix
$$
for nontrivial linear functionals $w_1,\dots,w_r$ on $\t$.
Note that the subspaces $V_i$ for $i>0$ are not
uniquely determined in general, and that the linear functionals
$\om_i$ are determined only up to sign.

By considering geodesics through $m$ (see \cite{Ko}), it
can be shown that the connected component of the fixed
point set of $T$ containing $m$ is a submanifold of $M$ --- the
geodesics through $m$ in the direction of $V_0$ give a local chart
(via the exponential map). This argument also shows that the
submanifold is totally geodesic, with tangent space $V_0$ at $m$.

Up to this point, we have used only the fact that 
the action of $T$ preserves the metric. Since $T$ preserves
the complex structure, 
$J$ commutes with the linear transformations $\Th_m(t)$ 
and $\th_m(t)$. We may therefore choose the decomposition 
so that $J$ preserves
each subspace $V_i$.  It follows that
each connected component of the fixed
point set of $T$ is actually a K\"ahler manifold, and
that (for $i>0$) we may write
$$
J\vert_{V_i} = \pm
\pmatrix
0 & 1 \\
-1 & 0
\endpmatrix
$$
An additional consequence of the existence of $J$ is that each
$V_i$ acquires a natural orientation, and so the linear
functionals $\om_i$ are now determined uniquely (for a given choice
of subspaces $V_i$.

We now turn to the computation of the Hessian $H$ of $f^X$.
In section 1.2 we gave a definition in terms of local coordinates,
but this is equivalent to the following more invariant definition
(see page 4 of \cite{Mi1}). For any vectors $V,W\in T_mM$, 
we have
$$
H(V,W)=\tilde V(\tilde W(f^X))(m)= \tilde V df^X(\tilde W)(m)
$$
where $\tilde V, \tilde W$ are any extensions of $V,W$ to
local vector fields on $M$. We therefore have
$H(V,W) = \tilde V \om(X^\ast,\tilde W)(m)
= \tilde V \al(X^\ast) (m)$, where $\al = -i_{\tilde W}\om$
(and $i$ denotes interior product). With a suitable choice
of the extensions $\tilde V, \tilde W$, the well known formula for
$d\al$ shows that $H(V,W)=-\om(\tilde W,[\tilde V,X^\ast])$.

From the definition of $\th$, it follows that 
$\th(X)V=[\tilde V,X^\ast]_m$. Hence we obtain the formula
$$
H(V,W)=\lan W,-J \th(X) V\ran.
$$
On each $V_i$ with $i>0$, it follows that $H$ is equal
to the inner product times the nontrivial linear functional $\om_i$.
Hence the Hessian is nondegenerate on a space complementary
to $V_0$, i.e. $f^X$ is a Morse-Bott function.

The index of the Hessian at a critical point $m$ is even,
being twice the number of $\om_i$'s such that $\om_i(m)<0$.
This implies that the Morse-Bott function $f^X$ is perfect,
and so the proof of Theorem 4.1.1 is complete.
\qed\enddemo

The linear functionals $w_1,\dots,w_r$ on $\t$ which
appear in this proof are called the
(nonzero) weights of the action of $T$ at the
fixed point $m$.  Theorem 4.1.1 generalizes to 
the case of a symplectic manifold, as
was pointed out in \cite{Fr1}.  We shall not need this extra generality,
and in fact we shall make essential use of the complex structure $J$, so we
shall only consider K\"ahler manifolds here.

The formula $-\nabla f^X = J X^\ast$ leads to a geometrical
description of the flow lines of $-\nabla f^X$, as we shall explain next.

\proclaim{Lemma 4.1.2} The action of $T=S^1\times \dots\times S^1$ 
on $M$ extends to an action of the
complexified torus $\TC=\C^\ast\times\dots\times\C^\ast$ on $M$. The vector
$\i X\in \tc$ generates the vector field $J X^\ast$, i.e. $(\i X)^\ast=J X^\ast$.
The flow line $\ga$ of $-\nabla f^X$ passing through $m\in M$ is given by
the action of the subgroup $\i \R X$, i.e. $\ga(t)=\exp\,\i tX\cdot m$.
\endproclaim

\demo{Sketch of the proof}
The extension of the action is guaranteed by the fact that $T$
preserves the complex structure of $M$ (see \cite{At2}). A
direct construction of the extension may be given by using the
integral curves of $J X^\ast$. The fact that  $(\i X)^\ast=J X^\ast$
follows (tautologically) 
from this, as does the
required description of the flow line $\ga$.
\qed\enddemo

\no Thus, the flow lines of the gradient vector field associated
to the action of $T$ are given by the action
of a subgroup of the larger group $\TC$.  
{\it Conversely, whenever an \ll algebraic torus\rr 
$\C^\ast \times\dots\times \C^\ast$ acts complex analytically
on a compact K\"ahler manifold
then we obtain both a Morse-Bott function and its gradient flow
in the above manner.}  Many manifolds do admit such actions, for
example the Grassmannian $\grkcn$ (which is acted upon
naturally by the group of diagonal matrices with nonzero complex
diagonal entries).  This provides a simple explanation
of the rather tricky calculation of Lemma 3.1.2, by means of which
we identified the gradient flow lines for the height functions.
It suffices to assume that the 
action of $\C^\ast \times\dots\times \C^\ast$
preserves the complex structure, because an 
$S^1 \times\dots\times S^1$-invariant  K\"ahler metric can be
obtained by averaging the given  K\"ahler metric over the
(compact) group $S^1 \times\dots\times S^1$.

We shall now address a question which lies at the heart of
Morse theory:  {\it how are the gradient flow lines of a
Morse function $f:M\to\R$
arranged within $M$?}  For example, which pairs of critical points
are connected by a flow line? And by how many flow lines? These
questions can answered by explicit computation of stable and
unstable manifolds in the case of a height function on a
Grassmannian, as we did in \S 3, but in general no such
computation will be possible.  Our main theme from now on will
be to consider this question for Morse-Bott functions $f^X$ associated
to torus actions. 

From the lemma it follows that the behaviour of the 
flow lines is related to the geometrical
properties of the various orbits of the group $\TC$.  If
$M$ is an algebraic K\"ahler manifold, and the
action of $\TC$ is algebraic (as we shall assume from now on),
the closures
of such orbits are special algebraic varieties called toric
varieties.\footnote{This is essentially the definition of a
toric variety.} 
In general toric varieties have singularities, but they are
particularly amenable to study (see \cite{Od}, \cite{Fu2})
because they may be  characterized by purely
combinatorial objects, called \ll fans\rrr. Now, in the
K\"ahler situation at least, the fan is equivalent to a
more familiar combinatorial object, namely a convex
polyhedron in Euclidean space. The next theorem describes this polyhedron.

Before stating the theorem, we need to introduce the moment map
associated to the action of $T$ on $M$. This is the map
$$
\mu:M\to\t^\ast
$$
which is determined up to an additive constant
by the condition 
$$
d\mu(\quad )(Y)=\om(Y^\ast,\quad )
$$
for all $Y\in\t$. The definition of
this map comes from symplectic geometry and classical
mechanics, but it has a straightforward Morse-theoretic interpretation:
$$
\align
&\mu(\quad)(Y) \ \text{is the Morse-Bott function $f^Y$ associated to}\\
&\text{the sub-torus of $T$ which is generated by $Y$.}  
\endalign
$$
\no The fact that
our Morse-Bott function $f^X$ is not alone, but is accompanied by a whole
family of Morse-Bott functions $f^Y$ parametrized by $Y\in\t$, is
significant.  The moment map $\mu$ assembles these Morse-Bott
functions into a single vector-valued function.

\proclaim{Theorem 4.1.3 (\cite{At2},\cite{Gu-St})}
Let $O_m$ denote the closure in $M$ of the orbit of $m$ under $\TC$,
i.e. $O_m=\overline{\TC\cdot m}$. Then

\no (i) $\mu(O_m)$ is the convex hull of the finite
set $\{\mu(m)\in\t^\ast \st m \ \text{is a critical point of}\ f^X\}$,

\no(ii) the inverse image (under $\mu$) of each open face of
$\mu(O_m)$ is a single $\TC$-orbit in $O_m$, and

\no(iii) $\mu$ induces a homeomorphism $O_m/T \to \mu(O_m)$
(although the action of $T$ on $O_m$ is not necessarily free).
\endproclaim

\demo{Proof} See Theorem 2 of \cite{At2}.
\qed\enddemo

The simplest nontrivial example of this theorem is given by
the action of $\C^\ast$ on $\C P^1 (\cong S^2)$ by
$u\cdot [z_0;z_1]=[uz_0;z_1]$. The corresponding Morse-Bott
function $f^X:S^2\to\R$ is a height function, and there
are two isolated critical points: the maximum point
and the minimum point.  We
have $f^X=\mu$ in this situation, and $f^X(S^2)$ is obviously
the line segment joining the maximum and minimum values.

We are specifically interested in the stable and unstable
manifolds of $f^X$, and their intersections.  Theorem 4.1.3
leads to the following information about these spaces.  For
simplicity we shall assume that $M$ is actually a smooth projective
variety, i.e. an algebraic submanifold of some complex projective
space, with the induced K\"ahler structure.  Furthermore we
assume that $\TC$ acts on $M$ by projective transformations.  From the
method of section 3.5 we can then deduce that the closures of the
stable  and unstable manifolds of our Morse functions are irreducible
algebraic subvarieties.  

\proclaim{Theorem 4.1.4 (\cite{Li})}
Assume that $M$ is a smooth projective variety,
and that $\TC$ acts on $M$ by projective transformations.
Let $V$ be an irreducible algebraic subvariety of 
$M$ which is preserved by
the action of $\TC$.
Then $\mu(V)$ is the convex hull of the finite
set $\{\mu(m)\in\t^\ast \st  
\ \text{$m$ is a critical point of $f^X$ in $V$} \}$.
\endproclaim

\demo{Sketch of the proof} Let $P_1,\dots,P_s$ be the distinct images
under $\mu$ of the critical points of $f^X$ which lie in $V$.
(This is necessarily a finite set, as $f^X$ is constant on any
connected critical submanifold.) For any one-dimensional subalgebra
$\R Y$ of $\t$, the image of the continuous function $f^Y\vert_V$
is a closed finite interval in $\R$ (as $V$ is connected).  Moreover,
since $\TC$ preserves $V$, the ends of this interval (i.e. the maximum
and minimum values of $f^Y$) are of the
form $f^Y(P_i),f^Y(P_j)$, for some $i,j$.  But
$f^Y = \pi_Y \circ \mu$, where $\pi_Y:\t^\ast\to\R$ is given by evaluation
at $Y$. It follows that $\mu(V)$ is contained in the
convex hull of $P_1,\dots,P_s$.

Conversely, we must show that any point of this polyhedron
is contained in $\mu(V)$. Let $P_{i_1},\dots,P_{i_t}$ denote
the \ll exterior\rr points of the polyhedron.
For each $i_j$, choose  some $Y_j\in\t$ 
such that the function $f^{Y_j}\vert_{V}$
has $P_{i_j}$ as its absolute minimum value, occurring on a
critical set $V_j$, where $V_j=V\cap M_j$ for some connected
component $M_j$ of the fixed point set of $\TC$ on $M$.
(This may be done by choosing a \ll generic\rr $Y_j$ such that 
$P_{i_j}$ is the closest point
of the polyhedron to the linear functional $\lan Y_j,\ \ran$,
where $\lan\ ,\ \ran$ is an invariant inner product on $\t$.) 

Let $S^V_j=V\cap S_j$, where $S_j$ is the stable manifold of
$M_j$ (for $f^{Y_j}$).  Since $\TC$ preserves $V$,
$S^V_j$ is the stable \ll manifold\rr of $V_j$ (for $f^{Y_j}\vert_{V}$).
Now, the stable manifold decomposition of $M$ induces a decomposition of $V$.  
The closure (in $V$) of each piece of this decomposition
is a subvariety of $V$,
and precisely one such closure must be equal to $V$ since
$V$ is irreducible.  Denote this piece by $V\cap S$, where $S$
is the corresponding stable manifold in $M$ (for $f^{Y_j}$).

We claim that $S=S_j$.  (This would be obvious if $V$ were
a smooth subvariety.) Assume that $S$ and $S_j$ are not equal,
so that $S\cap S_j=\es$. Then the closure of $V\cap S$ in $V$
is disjoint from $V\cap S_j$, since the closure of $S$ in
$M$ is disjoint from $S_j$. However this contradicts the defining
property of $V\cap S$. 

It follows that the complement of $S_j$ in $V$ is a subvariety
of positive codimension.
The complement of the intersection of all such $S_j$ (for
$j=1,\dots,t$) is therefore also a subvariety of positive
codimension; in particular this intersection
is nonempty. For any point $v$ of the intersection, we have
$P_{i_1},\dots,P_{i_t} \in \overline{\TC\cdot v}$, because
of the description of the flow lines as orbits of subgroups of $\TC$.
Hence $\mu(V)$, which contains $\mu(\overline{\TC\cdot v})$, must contain
the convex hull of $P_1,\dots,P_s$, by Theorem 4.1.3.
\qed\enddemo

\no The case $V=M$ was stated and proved as one of the main
theorems of \cite{At2} and \cite{Gu-St}. Although the proofs given
there were direct, the possibility of deducing the result from
Theorem 4.1.3 was in fact mentioned in \cite{At2}.   Various 
special cases of this result were already known --- we
refer to \cite{At2} and \cite{Gu-St} for further information.

Theorem 4.1.4 suggests the possibility that the behaviour of
the gradient flow lines of $f$ may be encoded by combinatorial
information within the polyhedron $\mu(M)$.  
We shall investigate this phenomenon, 
starting with the familiar case
of height functions on Grassmannians from \S 3. Initial
work in this direction was done in \cite{Ge-Ma}, where
some of the results of \cite{At2}, \cite{Gu-St} were anticipated
in the case of a Grassmannian.

Let $M=\grkcn$ and let $T_n$ be the group of diagonal
$n\times n$ matrices whose diagonal entries are complex
numbers of unit length, and 
let $T_n^{\C}$ be its complexification, i.e. 
the group of diagonal
$n\times n$ matrices with nonzero complex
diagonal entries.
We have a natural action of $T_n$ (or $T_n^{\C}$)
on $M$.
The vector $X=\i(a_1,\dots,a_n)$ in the Lie algebra $\t_n$ is
a generator of $T_n$ if and only if $\R(a_1,\dots,a_n)\cap \Z^{n}=\{0\}$,
i.e. the line through $(a_1,\dots,a_n)$ has \ll irrational slope\rrr.

From the above general theory, any choice of $(a_1,\dots,a_n)$ gives
rise to (1) an action of $T$ (the sub-torus of $T_n$ generated
by $Y=\i(a_1,\dots,a_n)$), and
(2) a Morse-Bott function $f^Y$, whose critical points are the
fixed points of $T$. From Lemma 3.1.2 we see that this
function is (up to an additive constant) the height function
on $\grkcn$ defined in \S 3.
The fixed points of $T$ are of course the $k$-planes which
can be spanned by eigenvalues of $\diag(a_1,\dots,a_n)$.

From the formula for $f^X$ in section 3.1, it follows that the moment map
is given explicitly by 
$$
\mu(V) = \text{diagonal part of}\ \pi_V.
$$
(We identify $\t^\ast$ with $\t\cong \R^n$ by using the 
(restriction of the) standard
Hermitian product on $\C^n$.)
In particular,
$$
\mu(V_{\u})=e_{\u},
$$
where $e_{\u} = e_{u_1} + \dots + e_{u_k}$ (and $e_{i}$ denotes the
$i$-th basis vector of $\C^n$). Since every critical manifold
of $f^X$ contains at least one point of the form $V_{\u}$,
the images of the critical points of $f^X$ under $\mu$ are
precisely the points $e_{\u}$. Hence the general theory gives
$$
\mu(\grkcn)= \ \text{convex hull of}\ 
\{e_{\u}= e_{u_1} + \dots + e_{u_k} \st 1\le u_1<\dots< u_k \le n \}.
$$

\proclaim{Example 4.1.5} 
\rm
Let $k=1$, so that $\grkcn=\C P^{n-1}$. In this case the
formula for $f^X$ in section 3.1 
shows that the moment map is given even more
explicitly by
$$
\mu([z])=(\vert z_1\vert^2,\dots,\vert z_n\vert^2) /
\sum_{i=0}^n\vert z_i\vert^2.
$$
It follows immediately from this formula that $\mu(\C P^{n-1})$
is the convex hull of the basis vectors $e_1,\dots,e_n$.
Since the (closures of the) 
stable and unstable manifolds of $f^X$ are respectively
of the form $[*;\dots;*;0;\dots;0]$ and $[0;\dots;0;*;\dots;*]$
(see section 3.1), their images under $\mu$ are also clear.
For example, in the diagram below, we illustrate the images under $\mu$
of the closures of the stable manifolds of $V_1$, $V_2$, $V_3$,
for a Morse function $f^X:\C P^2\to \R$.

$${}$$
$${}$$
$${}$$
$${}$$
$${}$$
$${}$$
$${}$$
$${}$$
$${}$$

\no In fact, for $\C P^{n-1}$, all assertions of Theorems 4.1.3
and 4.1.4 may be verified directly, from the formula for $\mu$.
\qed\endproclaim

\proclaim{Example 4.1.6}
\rm
Let us consider what the general theory says in the case of
$Gr_2(\C^4)$. We have already investigated this space in some
detail in \S 3 (starting with Example 3.4.1). The image of
$\mu$ is the convex hull of the six points $e_i + e_j$ in
$\R^4$ (with $1\le i<j \le 4$). It may be verified that
this polyhedron is a regular octahedron.

$${}$$
$${}$$
$${}$$
$${}$$
$${}$$
$${}$$
$${}$$
$${}$$
$${}$$

\no The heavy lines represent the partial order shown
in the diagram following Definition 3.4.2. Combining the
calculations of the stable manifolds in \S 3 with the
statement of Theorem
4.1.4, we see that the image under $\mu$ of the closure of the
stable manifold of a critical point $V_{\u}$ is
the convex hull of those vertices which are greater than
or equal to $e_{\u}$ in this partial order. 

With a little more work, it is possible to verify the predictions
of Theorem 4.1.3 in this case. First it is necessary
to identify the various possible types of the closures of the orbits of
the group $(\C^\ast)^4$, to which Theorem 4.1.3 will apply. 
Zero-dimensional
orbits, i.e. points, correspond to the vertices of the
octahedron. One-dimensional orbits (necessarily isomorphic
to $\C P^1$) correspond to the edges. Two-dimensional orbits
are of two types: copies of $\C P^2$, which correspond to faces,
and copies of $\C P^1\times \C P^1$, which correspond to
squares spanned by sets of four coplanar vertices.
Finally, there are two types of three-dimensional orbits,
whose images under $\mu$ are either half of the octahedron
or the entire octahedron. The first type is represented by
the Schubert variety $X_{(2,4)}$ (see section 3.4), and the
second by the famous \ll tetrahedral complex\rr (see \cite{Ge-Ma}).
\qed\endproclaim 

Returning to the general case of $\grkcn$, we mention that
the partial order on the critical points of $f^X:\grkcn\to\R$
(or equivalently, on the vertices of the polyhedron)
may be specified purely algebraically, in terms of the action
of the symmetric group $\Sigma_n$ (the Weyl group of $\glnc$).
This is explained in \cite{At2}.  

In \cite{Ge-Se} and \cite{Ge-Go-Ma-Se}, a detailed
study is made of the various sub-polyhedra of $\mu(\grkcn)$
which arise from Schubert varieties associated to the Morse
functions $f^X$ and their intersections. These are characterized
in terms of the combinatorial concept of a matroid.

By representing the stable and unstable manifolds of the
function $f^X$ as sub-polyhedra of the polyhedron $\mu(M)$,
we have in principle solved the problem of understanding 
the behaviour of the flow lines of $-\nabla f^X$. The
practical value of this solution depends on being able
to extract this information in an efficient manner, and
this may not always be possible. 
However, there are two general situations
where reasonably explicit results may be expected, namely
for generalized flag manifolds and for toric manifolds. The
essential phenomenon in both these situations is that
{\it there exists an orbit of the complex algebraic torus $\TC$
whose closure contains all critical points of the Morse functions
associated to the torus action.} We shall call
this kind of $\TC$-action a {\it complete torus action.}

\proclaim{Example 4.1.7}
\rm
A generalized flag manifold is by definition the
quotient of a complex semisimple Lie group $\GC$ by a
parabolic subgroup $P$. 
Height functions on generalized flag manifolds --- generalizing
the Morse functions on Grassmannians in \S 3 --- were
first studied by Bott (see \cite{Bo2} and the 
article of Bott in \cite{At1}). Such functions
are associated to torus actions (namely a maximal
torus of $\GC$) in the manner described at the
beginning of this section.  The index and, in the case
of a Morse-Bott function, the nullity of a critical
point may be computed in terms of the weights of the
torus action (i.e. in terms of the roots of $\GC$). 
Each stable or unstable manifold is an orbit of
a certain subgroup of $\GC$, as in the case $\GC/P=\grkcn$.
Indeed the decomposition of
$\GC/P$ into stable (or unstable) manifolds coincides
with the well known \ll Bruhat decomposition\rr of $\GC/P$, 
a fact which was proved in \cite{Pk}
as well as in various later papers
(e.g. \cite{At2}). A brief summary of this theory may be found
in the Appendix of \cite{Gu-Oh}.

The image of the moment map for the action of a
maximal torus on a generalized flag manifold was worked out
in \cite{At2}, generalizing earlier work of Kostant.
The polyhedron can be described (see \cite{Ge-Se})
as the convex hull of the {\it weights} of an
irreducible representation of $G$; it is well known that
$\GC/P$ is the projectivized  orbit of the maximal
weight vector of a suitable representation (the generalized
Pl\"ucker embedding).
Schubert varieties in generalized flag manifolds
have been extensively studied from the point of view
of Lie theory and algebraic geometry (see \cite{Hl}
for an introduction and further references).
The sub-polyhedra obtained by taking the images of
(various intersections of)
Schubert varieties in generalized flag manifolds have
been characterized in combinatorial terms in \cite{Ge-Se},
generalizing the results mentioned earlier for Grassmannians.
The homology classes represented by Schubert varieties have
also been investigated thoroughly. In \cite{Be-Ge-Ge} 
these homology classes are related to the well known description
of the cohomology ring due to Borel, by making use of the
generalized  Pl\"ucker embedding.
A brief explanation of the latter work can be found in \cite{Se}.
\qed\endproclaim 

\proclaim{Example 4.1.8}
\rm
From the general theory of toric varieties,
it is well known that a (smooth) toric variety with
a K\"ahler metric is
entirely determined by the image of its moment mapping.
In particular, the behaviour of the flow lines of
a Morse-Bott function associated to the given
torus action is represented faithfully in the
momentum polyhedron.  Various geometrical and
topological invariants of such manifolds have been
computed explicitly in terms of this polyhedron;
details can be found in \cite{Od} and \cite{Fu2}
(see also \cite{Au} and \cite{De}).
\qed\endproclaim 

There is some intersection between Examples 4.1.7
and 4.1.8, as one may consider the toric varieties obtained
as the closures of the orbits of a maximal torus of
$\GC$ acting on $\GC/P$. These are singular varieties,
in general; they have been studied in \cite{Ge-Se} and later
in \cite{Fl-Ha}, \cite{Da}, \cite{Ca-Ku}.

We shall now change our point of view slightly by focusing on
the family of Morse-Bott functions $f^Y$ (parametrized by
$Y\in\t$), rather than the single Morse function $f^X$ (corresponding
to a generator $X$ of $\t$). This theme will reach maturity
in section 4.3, so as motivation for this we shall consider again the
problem of computing the cohomology ring $H^\ast(Gr_2(\C^4))$
(cf. section 3.6).

Theorem 4.1.4 gives a representation of the (image under $\mu$ of
the) Schubert variety $V=X^Y_{\u}$ for the Morse function
$f^Y:\grkcn\to\R$. 
Namely, the image under $\mu$ of $V$ is the convex
hull of those points $e_{\u}$ such that
$V_{\u}\in V$.
This representation also applies to the
(irreducible components of the) 
variety $V=X^{Y_1}_{\u_1}\cap\dots\cap X^{Y_k}_{\u_k}$.
In section 3.6 we saw that the problem of
calculating products in cohomology can in principle be
reduced to the problem of calculating 
intersections of (pairs or) triples
of Schubert varieties 
$X^{Y_1}_{\u_1}\cap X^{Y_2}_{\u_2}\cap X^{Y_3}_{\u_3}$.
\ll In principle\rr means\footnote{It can be 
shown that all intersections are indeed transverse,
by using the fact that the stable manifolds are orbits of certain
subgroups of $\glnc$; see \cite{Kt}. This is not quite enough,
for we must show in addition that there exist sufficiently many 
representatives of all cohomology classes. This will be obvious
in the examples we consider, however.}
\ll providing that all necessary
triple intersections are transverse\rrr.
To be more precise, we need to find all zero-dimensional triple
intersections. By Theorem 4.1.4 (or by the Pl\"ucker embedding
argument of section 3.6), such an intersection either consists
of a single point or is empty. To determine which is the case,
we just need to know which cohomology classes are represented
by which sub-polyhedra of $\mu(\grkcn)$. Let us consider two
examples, $\C P^2$ and $Gr_2(\C^4)$.

\proclaim{Example 4.1.9}
\rm
The cohomology ring $H^\ast(\C P^2)$ has additive generators
in dimensions $0$, $2$, and $4$; let us denote these
respectively by $1$, $A$, and $B$. They are dual to the
fundamental homology classes of the Schubert varieties for
a fixed Morse function $f^Y$. By varying $Y$ in $\t$,
we arrive at the following representation of these 
cohomology classes on the triangular region $\mu(\C P^2)$:

$${}$$
$${}$$
$${}$$
$${}$$
$${}$$
$${}$$
$${}$$
$${}$$
$${}$$

The structure of the cohomology ring $H^\ast(\C P^2)$ is
determined completely by the product $A^2$, and from the diagram
the intersection number of the dual Schubert varieties is $1$.
Hence $A^2=1B=B$, as expected.
\qed\endproclaim 

\proclaim{Example 4.1.10}
\rm
The cohomology ring $H^\ast(Gr_2(\C^4))$ has the six additive generators
described in section 3.6.  These are represented on the 
octahedron of Example 4.1.6 as follows:
$$
\align
z_{(3,4)} &\in H^0:
\quad\text{the octahedron}\\
z_{(2,4)} &\in H^2:
\quad\text{the half octahedra}\\
z_{(1,4)}, z_{(2,3)} &\in H^4:
\quad\text{alternate faces}\\
z_{(1,3)} &\in H^6:
\quad\text{the edges}\\
z_{(1,2)} &\in H^8:
\quad\text{the vertices}
\endalign
$$

In the diagram below, the (four) faces which represent the cohomology class
$z_{(2,3)}$ are shaded.

$${}$$
$${}$$
$${}$$
$${}$$
$${}$$
$${}$$
$${}$$
$${}$$
$${}$$

It is now a simple matter to read off all zero-dimensional
double and triple intersections. For example, 
$z_{(1,4)}z_{(2,4)}z_{(2,4)}$ is represented by the
intersection of two half octahedra and a face.
Any such intersection {\it giving a finite number of points}
gives precisely one point. So the product is equal to the
generator of $H^8$ --- as we found by a much more laborious calculation in
Example 3.6.4.  By contemplating the above diagram we can determine
the entire cohomology ring of $Gr_2(\C^4)$!
\qed\endproclaim 

To end this section, we emphasize two advantages of 
having a Morse-Bott function which is associated to a torus action.
First, having a group action has computational advantages, as we
have seen in the description of the gradient flow and the 
identification of the index (and nullity) of a critical point.  
Second, the torus action gives more than a single Morse-Bott function;
it gives a whole family of related Morse-Bott functions, and this
family can give more information than one of its members
(as in the above calculation of the cohomology ring).

\noindent{\it Additional comments (May, 2000):} The theory of
equivariant cohomology provides an algebraic explanation
for the success of the above method of computing $H^\ast(M)$ from
a (complete) torus action (see M. Goresky, R. Kottwitz, and R. MacPherson,
{\it Equivariant cohomology, Koszul duality,
and the localization theorem,} Invent. Math. 131 (1998), 25-83).
The equivariant cohomology ring $H^\ast_T(M)$ (see section 2.7)
is in this case isomorphic as a module over $H^\ast(BT)$ to the
tensor product $H^\ast(M)\otimes H^\ast(BT)$. The localization
theorem for equivariant cohomology expresses 
products of equivariant cohomology classes in terms of the fixed
point data of the torus action (this is another manifestation of
Morse theory). As a result, it is possible to describe $H^\ast(M)$
explicitly as a quotient of a polynomial ring, in terms of
this data.

\subheading{4.2 The Witten complex}

In view of the importance of the gradient flow lines 
of a Morse function, it is perhaps not surprising that
the basic theorems of Morse theory may be developed
entirely from this point of view.  In fact,
much stronger results (than the traditional \ll Morse inequalities\rrr)
are possible, as we shall see in this section and the next one.

In this section we consider the goal of computing the
homology groups of a manifold $M$. Traditionally, this is
possible only for a {\it perfect} Morse function.  
However, if we assume that
$f:M\to\R$ is a Morse-Smale function (i.e. the stable and
unstable manifolds of $f$ intersect transversely, as in
Definition 1.4.5), then the homology may be calculated
whether $f$ is perfect or not.  This method became widely
understood only in the 1980's, through the work of Witten 
and Floer (see \cite{Wi} and \cite{Fl}).
It is easy to describe: one constructs a certain chain
complex of abelian groups, the \ll Witten complex\rrr, 
whose homology groups turn out to be the homology groups of $M$.

Let us assume first that $M$ is oriented.
The abelian
groups are defined in terms of the critical points of the Morse-Smale
function $f$ by
$$
C_i = \ \text{free abelian group on the set of
critical points of index $i$}.
$$
(Since $f$ is a Morse function and $M$ is compact,
the groups $C_i$ have finite rank.)
The boundary operators $\b_i : C_i \to C_{i-1}$
are defined in terms of the gradient flow lines of
$-\nabla f$ as follows. Let $m$ be a critical point
of $f$ of index $i$ (i.e. a generator of $C_i$). Then
$$
\b_i(m) = \sum_{\ga} e(\ga) m_{\ga},
$$
where 

\no(1) the sum is over all flow lines $\ga$ such that
$$
\lim_{t\to -\infty} \ga(t) = m, \quad
\lim_{t\to \infty} \ga(t) = \ \text{(a critical point) $m_{\ga}$
of index}\ i-1, 
$$

\no(2) $e(\ga)$ is either $1$ or $-1$, the choice depending on
whether $\ga$ \ll preserves or reverses orientation\rrr.

Some explanation of (1) and (2) is necessary.  First, since
the Morse-Smale condition gives
$$
\dim F(m,m_{\ga}) = 1,
$$
it follows that there are only finitely many such $\ga$,
hence the sum is finite. 
Second, to define the sign of $e(\ga)$, we first choose arbitrary
orientations of the unstable manifolds. Since $M$ is oriented,
and since the stable and unstable manifolds intersect
transversely, we may then assign orientations to the stable manifolds
in a consistent manner. The manifold $F(m,m_{\ga})$ itself
thus acquires an orientation. We define $e(\ga)$ to
be $1$ if the natural orientation of $\ga$ agrees with
its orientation as a component of $F(m,m_{\ga})$; otherwise
we define $e(\ga)$ to be $-1$.

It can be shown that $(C_{\ast},\b_{\ast})$ is a chain complex,
i.e. that $\b_{i-1}\circ \b_i=0$ for all $i$.  
For this, and for the proof of the
next theorem, we refer to section 2 of \cite{Au-Br}, where
a detailed discussion can be found.

\proclaim{Theorem 4.2.1} The homology groups of the
chain complex $(C_{\ast},\b_{\ast})$ are isomorphic to the homology
groups of $M$.
\qed\endproclaim

If $M$ is orientable, as in the above definition, then
homology groups with coefficients in $\Z$ are obtained. 
If $M$ is not orientable, then some modifications to
the definition are necessary.  The simplest way to
do this is to work over $\Z/2\Z$ instead of $\Z$,
i.e. to replace $C_i$ by $C_i\otimes \Z/2\Z$ and
then to define $e(\ga)=1$ for all $\ga$. In this case
the theorem is true for homology groups with coefficients in 
$\Z/2\Z$.

It is possible to formulate and prove a similar theorem
for the cohomology groups of $M$,
using the dual chain complex (see \cite{Au-Br}). 

The Morse inequalities (Theorem 2.3.1) follow from the
statement of the
above theorem by a purely algebraic argument (using
$\rank H_i(M) = \rank \Ker \b_i - \rank \Im \b_{i+1}$
and
$\rank C_i = \rank \Ker \b_i + \rank \Im \b_{i}$).
The \ll lacunary principle\rrr, that $f$ is necessarily
perfect if all its critical points have even index, also
follows immediately, since in this case we have $\b_{i} = 0$
for all $i$.

\proclaim{Example 4.2.2}
\rm
Consider the Morse function on the circle with three
local maxima and three local minima depicted in Example 2.2.2.

$${}$$
$${}$$
$${}$$
$${}$$
$${}$$
$${}$$
$${}$$
$${}$$
$${}$$

This is of course not a perfect Morse function. But the
homology groups may be calculated by using the
Witten complex. If we choose the \ll clockwise\rr
orientation on $S^1$ and on all stable manifolds,
then the map $\b_1:C_1\to C_0$ is given by:
$$
\b_1 F = -A + C, \quad \b_1 D = A - B, \quad \b_1 E=  B-C.
$$
The kernel and cokernel of $\b_1$ are therefore both
isomorphic to $\Z$.
\qed\endproclaim

\proclaim{Example 4.2.3}
\rm
Consider the Morse function $f:\R P^n\to\R$
defined by 
$$
f([x_0,\dots,x_n])=\sum_{i=0}^n c_i\vert x_i\vert^2 /
\sum_{i=0}^n\vert x_i\vert^2  
$$
with $c_0>\dots> c_n$, which we met in Example 2.3.6. 
Let us calculate the homology of $\R P^n$ by using
the Witten complex. First we take coefficients in $\Z/2\Z$,
to avoid the problem of dealing with orientations.

The critical points are the coordinate axes $V_0,\dots,V_n$,
and the indices are (respectively) $n,n-1,\dots,0$. 
Thus $C_i = \Z/2\Z\, V_{n-i}$ for $0\le i\le n$.
As in the case of $\C P^n$ in \S 3, we may identify the
stable and unstable manifolds explicitly. In particular,
we see that the space $F(V_i,V_{i+1})$ of points
on flow lines from $V_{i}$ to $V_{i+1}$ is of the form
$$
\{ [0;\dots;0;\ast;\ast;0;\dots;0] \in \R P^n
\st \ast \in \R \}.
$$
In other words, it is a copy of $\R P^1 \cong S^1$, and
so there are precisely two such flow lines. Thus, every
homomorphism $\b_i$ is zero, and our Morse function is
perfect.

With integer coefficients the situation is more complicated,
particularly since we have not defined the Witten complex
(over $\Z$) for a nonorientable manifold, and it is well known
that $\R P^n$ is orientable only when $n$ is odd. However,
the Witten complex can in fact be defined for any manifold
(see section 2.1 of \cite{Kt}), and it turns out that for $\R P^n$
the maps $\b_i$ are given by
$$
\b_i(n)=
\cases
2n \ \text{if $i$ is even}\\
0 \ \text{if $i$ is odd}
\endcases
$$
This gives the integral homology groups of $\R P^n$.
\qed\endproclaim

\proclaim{Example 4.2.4}
\rm
Any function associated to the standard torus action
on a generalized flag manifold (see Example 4.1.7) 
has trivial Witten complex, since
all critical points have even index. However,
the real analogues of these complex manifolds (which
include $\R P^n$ and more generally the real Grassmannians)
give rise to nontrivial Witten complexes, and it
may be expected that these complexes are determined by 
the same combinatorial information that describes the
behaviour of the flow lines. This is indeed the
case; full details may be found in \cite{Kt}.
Extensive earlier work on the Morse theory of these spaces
can be found in \cite{Bo-Sa}, \cite{Tk}, \cite{Tk-Ko},
\cite{Fr2}, \cite{Du-Ko-Va}. The convexity results of \cite{At2},
\cite{Gu-St} were extended to these spaces in \cite{Du}.
\qed\endproclaim

The above description of the Witten complex, phrased
in the language of differential topology, is closely
related to earlier work of Smale and Thom
(see \cite{Fs}, \cite{Mi2}, \cite{Sm}), in which
the groups $C_i$ appear as the relative homology
groups of the pair $(M_i,M_{i-1})$, for a suitable
filtration $\{ M_i \}$ of $M$, and the
maps $\b_i$ appear as the connecting homomorphisms
in the homology exact sequence. (The space $M_i$ is
obtained by taking all cells of dimension less than or
equal to $i$ in the usual Morse decomposition of $M$.)
Witten's original motivation was actually
quite different, as it arose from quantum theory.
A brief description of Witten's point of view, together with further
historical information, can be found in \cite{Bo4}.

Finally we mention that an extension of the theory
to the case of
Morse-Bott functions is given in \cite{Au-Br}.
Another reference where full details of the material of this section 
can be found is the book \cite{Sc}.

\subheading{4.3 Morse theory as a topological field theory}

In section 4.1 we have seen that it can be useful to study 
{\it families} of Morse functions on a given compact manifold
$M$. This can be taken as motivation for the \ll field-theoretic\rr 
approach to Morse theory of \cite{Be}, \cite{Be-Co}, \cite{Co-Jo-Se1},
\cite{Co-Jo-Se2},
\cite{Fy}, \cite{Fy-Se}, so called because it is based on
similar constructions in gauge theory.  As such constructions
tend to involve rather elaborate preparations, we shall just 
give an informal description here.

The basic ingredient is a certain \ll moduli space\rr $\Cal M(\Ga)$,
which is a device for counting configurations of flow lines.
The definition of this space depends on $M$ and on a choice of
an oriented connected graph $\Ga$. We assume that $\Ga$ has
$n_1$ edges parametrized by $(-\infty,0]$ (\ll incoming
edges\rrr), 
$n_2$ edges parametrized by $[0,1]$ (\ll internal
edges\rrr), and
$n_3$ edges parametrized by $[0,\infty)$ (\ll outgoing
edges\rrr). In the example below, we have
$n_1=2$, $n_2=1$, and $n_3=3$.

$${}$$
$${}$$ 
$${}$$
$${}$$ 
$${}$$
$${}$$ 
$${}$$
$${}$$ 

An element of $\Cal M(\Ga)$ is a \ll configuration of
flow lines of Morse-Smale functions on $M$, modelled on $\Ga$\rrr,
i.e. a continuous map $\Cal F:\Ga \to M$ such that, on each edge of $\Ga$,
$\Cal F$ is (part of) a flow line of a Morse-Smale function on $M$.
Thus, for the graph illustrated above, $\Cal F(\Ga)$ might look like
the diagram below, where the white dots are critical points approached
by the incoming and outgoing flow lines.

$${}$$
$${}$$ 
$${}$$
$${}$$ 
$${}$$
$${}$$ 
$${}$$
$${}$$ 

Let $f_1^{\Cal F}, \dots, f_{n_1+n_2+n_3}^{\Cal F}$
be the Morse-Smale functions in the definition of $\Cal F$,
and let $a_1^{\Cal F},\dots, a_{n_1}^{\Cal F},
a_{n_1+n_2+1}^{\Cal F}, \dots, a_{n_1+n_2+n_3}^{\Cal F}$
be the critical points which are approached by the incoming
and outgoing flow lines in that definition
(the white dots in the diagram).
For any $(n_1+n_2+n_3)$-tuple of functions
$g=(g_1,\dots,g_{n_1+n_2+n_3})$ on $M$, we define
$$
\Cal M_g(\Ga) =
\{ \Cal F \in \Cal M(\Ga) \st f_i^{\Cal F} = g_i \}.
$$
For any $(n_1+n_3)$-tuple of points 
$b=(b_1,\dots,b_{n_1},b_{n_1+n_3+1},\dots,b_{n_1+n_2+n_3})$ of $M$,
we define
$$
\Cal M_g(\Ga;b) =
\{ \Cal F \in \Cal M_g(\Ga) \st a_i^{\Cal F} = b_i \}.
$$
(These definitions are informal versions of the
precise definitions in \cite{Be-Co}.)

{\it We shall assume from now on that all stable and unstable manifolds of} 
$f_1^{\Cal F}, \dots, f_{n_1+n_2+n_3}^{\Cal F}$
{\it intersect transversely;} in particular all these
functions are Morse-Smale functions. Under this assumption,
it can be shown that

\no(i)  $\Cal M_g(\Ga;b)$ is a smooth manifold

\no(ii)  $\Cal M_g(\Ga;b)$ is oriented if $M$ is oriented, and

\no(iii)  $\dim \Cal M_g(\Ga;b) =
\sum_{i=1}^{n_1} \index b_i  - \sum_{i=1}^{n_3} \index b_{n_1+n_2+i}
- (\dim M)(\dim H_1(\Ga;\R) + n_1 - 1)$.

\proclaim{Example 4.3.1}
\rm
Let $\Ga$ be the graph below  with $n_1=n_3=1$, $n_2=0$:

$${}$$
$${}$$ 
$${}$$
$${}$$ 
$${}$$

\no Let $g=(g_1,g_2)$, $b=(b_1,b_2)$. Then the points of $\Cal M_g(\Ga;b)$
are in one-to-one correspondence with the points of 
$U_{b_1}^{g_1}\cap S_{b_2}^{g_2}$, where $U_{b_i}^{g_i}$ is
the unstable manifold of the critical point $b_i$ of $g_i$, and
$S_{b_i}^{g_i}$ is the stable manifold. From the transversality
assumption we have (see section 1.4)
$$
\codim U_{b_1}^{g_1}\cap S_{b_2}^{g_2} =
(\dim M - \index b_1) + (\index b_2).
$$
This checks with the general formula above, i.e.
$\dim \Cal M_g(\Ga;b) = \index b_1 - \index b_2$.
\qed\endproclaim

\proclaim{Example 4.3.2}
\rm
Let $\Ga$ be the graph below with $n_1=2$, $n_2=0$, and $n_3=1$:

$${}$$
$${}$$ 
$${}$$
$${}$$ 
$${}$$

\no Let $g=(g_1,g_2,g_3)$, $b=(b_1,b_2, b_3)$. In this
situation the points of $\Cal M_g(\Ga;b)$ correspond to
points of $U_{b_1}^{g_1}\cap U_{b_2}^{g_2}\cap S_{b_3}^{g_3}$.
Transversality implies that
$$
\codim U_{b_1}^{g_1}\cap U_{b_2}^{g_2}\cap S_{b_3}^{g_3} =
(\dim M - \index b_1) + (\dim M - \index b_2) + (\index b_3).
$$
Again this is consistent with the general formula.
\qed\endproclaim

We come now to the main part of the construction, which will
associate to each graph $\Ga$ a topological invariant of $M$.
This will make use of the Witten complex $(C_{\ast}(f), \b_{\ast})$
of a (Morse-Smale) function $f$, which was defined in the
last section.

For a fixed choice of $g$ (as above), we define
$$
q_g(\Ga)=
\sum_{\Cal F,b} e(\Cal F)\, b_1\otimes \dots\otimes b_{n_1}
\otimes b_{n_1+n_2+1} \otimes\dots\otimes b_{n_1+n_2+n_3},
$$
where the sum is over all $\Cal F, b$ such that
$\Cal F\in \Cal M_g(\Ga;b)$ and $\dim \Cal M_g(\Ga;b) = 0$.
If $M$ is oriented, so is the zero-dimensional manifold
$\Cal M_g(\Ga;b)$, and $e(\Cal F)$ is plus or minus one,
according to the orientation of $\Cal F$ as a point of $\Cal M_g(\Ga;b)$.

Thus, if $M$ is oriented, we may regard $q_g(\Ga)$
as an element of
$$
\left(\bigotimes_{i=1}^{n_1} C_{\ast}(g_i)\right) \ \otimes\ 
\left(\bigotimes_{i=n_1+n_2+1}^{n_1+n_2+n_3} C^{\ast}(g_i)\right),
$$
where $C^{\ast}(g_i)$ is the dual complex to $C_{\ast}(g_i)$.
If $M$ is not oriented, then the definition of $e(\Cal F)$ must 
be modified in the same way as the Witten complex.

The construction of $q_g(\Ga)$ depends on the choice
of the $(n_1+n_2+n_3)$-tuple of Morse-Smale functions
$g$, of course. Nevetheless, it turns out that this choice,
and all other choices necessary for a rigorous definition of
$q_g(\Ga)$, are irrelevant:

\proclaim{Theorem 4.3.3 (\cite{Be}, \cite{Be-Co}, \cite{Fy})}
The element $q_g(\Ga)$ is annihilated by the
(appropriate extension of) $\b_{\ast}$, and so we obtain a class
$$
[q_g(\Ga)] \in
\left(\bigotimes_{i=1}^{n_1} H_{\ast}(M)\right) \ \otimes\ 
\left(\bigotimes_{i=n_1+n_2+1}^{n_1+n_2+n_3} H^{\ast}(M)\right).
$$
This element depends only on the graph $\Ga$.
\qed\endproclaim

We denote the class $[q_g(\Ga)]$ by $q(\Ga)$.  We may regard
$q(\Ga)$ as an element of the ring
$\Hom(\bigoplus H^{\ast}(M),\bigoplus H^{\ast}(M))$,
i.e. as an {\it operation} on cohomology classes. Various
well known operations can be obtained this way, by
choosing suitable graphs.  We shall discuss the case of
the (cup) product operation in the following example.

\proclaim{Example 4.3.4}
\rm
We choose the graph illustrated below:

$${}$$
$${}$$ 
$${}$$
$${}$$ 
$${}$$

\no If $\Cal M_g(\Ga;b) 
(\cong U_{b_1}^{g_1}\cap U_{b_2}^{g_2}\cap U_{b_3}^{g_3})$
is nonempty and zero-dimensional, 
then by transversality we must have
$\sum_{i=1}^3 (\dim M - \index b_i)=\dim M$.
We obtain
$$
q(\Ga)\in
\bigoplus_{i_1+i_2+i_3=\dim M}
H_{\dim M - i_1}(M)\otimes H_{\dim M - i_2}(M)
\otimes H_{\dim M - i_3}(M),
$$
and hence also
$$
q(\Ga)\in
\bigoplus_{i_1+i_2+i_3=\dim M}
H^{i_1}(M)^\ast\otimes H^{i_2}(M)^\ast\otimes H^{i_3}(M)^\ast.
$$
From our calculations of the triple product operation
of $H^{\ast}(\grkcn)$ in section 3.6 and the end of
section 4.1, it is clear that
$q(\Ga)$ must be precisely that operation. This works
for any orientable manifold $M$, because the main ingredient
used in our calculation for $\grkcn$ was the existence of three
functions, all of whose stable and unstable manifolds
intersect transversely.  In the general case, however,
it is not easy to identify such \ll generic\rr functions
explicitly.  The special feature of $\grkcn$, and
indeed of any manifold with a torus action which satisfies
the condition stated before Example 4.1.7, is that one only
needs generic elements of the Lie algebra of the torus,
and the existence of these is guaranteed by the convexity theorem.
\qed\endproclaim

The product operation has been described by other authors in
terms of the Witten complex --- see \cite{Au-Br} and \cite{Vi}.

In all our examples so far we had $n_2=0$, and the
space $\Cal M_g(\Ga;b)$ was identified with a subspace
of $M$ itself. In general, $\Cal M_g(\Ga;b)$ may be
identified with a subspace of the $(n_2+1)$-fold product
$M\times \dots \times M$. 

The theory described in this section makes use only of those
$\Cal M_g(\Ga;b)$ which are zero-dimensional. In
\cite{Co-Jo-Se1},  \cite{Co-Jo-Se2} the higher dimensional
spaces are used to construct a much more complicated algebraic
object, which gives correspondingly more topological
information.

\subheading{4.4 Origins and other directions}

The study of critical points of functions on 
infinite-dimensional spaces (e.g. on function spaces --- the
Calculus of Variations) has been a guiding principle right from 
the beginning of Morse theory.  Rather surprisingly, perhaps,
the development of almost all of Morse theory 
has been prompted by infinite-dimensional examples! Since
the infinite-dimensional theory is much more complicated,
it is usually presented as a generalization
of traditional Morse theory, and we shall continue this
tradition by giving only a brief list of examples, almost as an
afterthought. Nevertheless, it is very likely that future
directions in Morse theory will be strongly influenced
by examples like these.

One of the earliest examples was the study of geodesics as
critical points of the length or energy functional on the
space of paths on a Riemannian manifold. Morse's idea of
using \ll broken geodesics\rr to reduce the problem to a
finite-dimensional problem is described in detail in \cite{Mi1}.
Subsequently, a general theory of Morse functions on Hilbert
manifolds was developed by Palais and Smale 
(see \cite{Pa1} and \cite{Pa2}), under the assumption of
\ll Condition (C)\rrr.  This condition, a substitute for
compactness, is satisfied in the case of the geodesic problem,
but unfortunately not in the case of many other important examples.
For example, it is not satisfied in the case of the energy
functional on the space of maps from a Riemann surface into
a Riemannian manifold. In this case the critical points are the
harmonic maps, which are closely related to minimal surfaces.

For the Yang-Mills functional on the space of connections
over a Riemann surface, a substitute for compactness was 
found by Atiyah and Bott (\cite{At-Bo}). This led to new
developments in finite-dimensional Morse theory, namely for the
functional $\vert \mu\vert^2$ where $\mu$ is the moment
map for the action of a (not necessarily abelian)
Lie group on a manifold $M$. In \cite{Ki}, Kirwan showed
that a version of Morse theory holds for this function,
even though it is not a Morse-Bott function.

A general approach to Morse theory for the Yang-Mills functional
(and other functionals, such as the energy functional on
maps of Riemann surfaces)
has been given by Taubes (see \cite{Tu} and also
the survey article \cite{Uh} of Uhlenbeck).

All these examples focus on the critical points of
a functional, but (as we have seen) it is possible to take a
different point of view by focusing on the flow lines.
Floer's idea of restricting attention to certain well-behaved
flow lines is such a case, and this was also 
motivated by an infinite-dimensional
example --- the \ll area functional\rr on closed paths in a
symplectic manifold $X$.  The flow lines in this case have
particular geometrical significance: they may be regarded as
\ll holomorphic curves\rr in $X$, where $X$ is given an almost 
complex structure compatible with its symplectic structure.
Such holomorphic curves arose in earlier work of Gromov, and 
the homology theory which is computed by the Witten
complex in this situation is known as Gromov-Floer theory.
(It is also closely related to quantum cohomology theory.)

Floer applied his theory to an apparently quite different
example, the Chern-Simons functional on the space of
connections on certain three-dimensional manifolds.  Again the
flow lines have a geometrical meaning: they are Yang-Mills instantons.
The homology theory arising here is called Floer homology.

\eightpoint \it
      
\no  Department of Mathematics, Graduate School of Science,
Tokyo Metropolitan University, Minami-Ohsawa 1-1, Hachioji-shi,
Tokyo 192-0397, Japan

\no martin\@comp.metro-u.ac.jp

\enddocument